\documentclass[11pt]{article}
\usepackage{amssymb}
\usepackage{fullpage,amsfonts,color,amsmath,amsthm,graphicx}

\newcommand{\pt}[1]{{\color{black} #1}}

\newcommand{\mm}[1]{{\color{black} #1}}
\newcommand{\mmm}[1]{{\color{black} #1}}
\newcommand{\mmmm}[1]{{\color{black} #1}}
\newcommand{\mmt}[1]{{\color{black} #1}}

\newcommand{\dm}[1]{{\color{black} #1}}
\newcommand{\dmm}[1]{{\color{black} #1}}
\newcommand{\dmmm}[1]{{\color{black} #1}}

\begin{document}

\newcommand{\N}{{\mathbb N}}

\def\a{\alpha}
\def\b{\beta}
\def\c{\chi}
\def\d{\delta}
\def\D{\Delta}
\def\e{\epsilon}
\def\f{\phi}
\def\F{\Phi}
\def\g{\gamma}
\def\G{\Gamma}
\def\k{\kappa}
\def\K{\Kappa}
\def\z{\zeta}
\def\th{\theta}
\def\Th{\Theta}
\def\la{\lambda}
\def\La{\Lambda}
\def\m{\mu}
\def\n{\nu}
\def\p{\pi}
\def\P{\Pi}
\def\r{\rho}
\def\R{\Rho}
\def\s{\sigma}
\def\S{\Sigma}
\def\t{\tau}
\def\om{\omega}
\def\Om{\Omega}
\def\smallo{{\rm o}}
\def\bigo{{\rm O}}
\def\to{\rightarrow}
\def\E{{\bf E}}
\def\ex{{\bf E}}
\def\cd{{\cal D}}
\def\rme{{\rm e}}
\def\hf{{1\over2}}
\def\R{{\bf  R}}
\def\cala{{\cal A}}
\def\cale{{\cal E}}
\def\Fscr{{\cal F}}
\def\cc{{\cal C}}
\def\calc{{\cal C}}
\def\calh{{\cal H}}
\def\call{{\cal L}}
\def\calr{{\cal R}}
\def\calb{{\cal B}}
\def\calw{{\cal W}}
\def\calz{{\cal Z}}
\def\bk{\backslash}

\def\out{{\rm Out}}
\def\temp{{\rm Temp}}
\def\overused{{\rm Overused}}
\def\big{{\rm Big}}
\def\notbig{{\rm Notbig}}
\def\moderate{{\rm Moderate}}
\def\swappable{{\rm Swappable}}
\def\candidate{{\rm Candidate}}
\def\bad{{\rm Bad}}
\def\crit{{\rm Crit}}
\def\col{{\rm Col}}
\def\dist{{\rm dist}}
\def\deg{{\rm deg}}
\def\whp{w.h.p.}
\def\cmin{c_{\rm min}}
\def\cmax{c_{\rm max}}

\newcommand\fix{{\rm FIX }}
\newcommand\fixx{{\rm FIX2 }}
\newcommand{\blank}{{\rm Blank }}

\newcommand{\Exp}{\mbox{\bf E}}
\newcommand{\var}{\mbox{\bf Var}}
\newcommand{\pr}{\mbox{\bf Pr}}
\newcommand{\Bin}{\mbox{\bf Bin}}
\newcommand{\Po}{\mbox{\bf Po}}

\newtheorem{lemma}{Lemma}
\newtheorem{theorem}[lemma]{Theorem}
\newtheorem{corollary}[lemma]{Corollary}
\newtheorem{claim}[lemma]{Claim}
\newtheorem{remark}[lemma]{Remark}
\newtheorem{observation}[lemma]{Observation}
\newtheorem{proposition}[lemma]{Proposition}
\newtheorem{definition}[lemma]{Definition}

\newcommand{\limninf}{\lim_{n \rightarrow \infty}}
\newcommand{\proofstart}{{\bf Proof.\hspace{2em}}}
\newcommand{\tset}{\mbox{$\cal T$}}
\newcommand{\proofend}{\hspace*{\fill}\mbox{$\Box$}\vspace{2ex}}
\newcommand{\bfm}[1]{\mbox{\boldmath $#1$}}
\newcommand{\reals}{\mbox{\bfm{R}}}
\newcommand{\expect}{\mbox{\bf Exp}}
\newcommand{\he}{\hat{\e}}
\newcommand{\card}[1]{\mbox{$|#1|$}}
\newcommand{\rup}[1]{\mbox{$\lceil{ #1}\rceil$}}
\newcommand{\rdn}[1]{\mbox{$\lfloor{ #1}\rfloor$}}
\newcommand{\ov}[1]{\mbox{$\overline{ #1}$}}
\newcommand{\inv}[1]{\frac{1}{ #1}}

\def\calc{{\cal C}}
\def\cald{{\cal D}}

\author{
Dieter Mitsche\thanks{Universite de Nice Sophia-Antipolis, Nice, France.}
\and
Michael Molloy\thanks{Department of Computer Science, University of Toronto, Toronto, ON, Canada. }
\and
Pawe\l{}~Pra\l{}at\thanks{Department of Mathematics, Ryerson University, Toronto, ON, Canada.}
}

\title{$k$-regular subgraphs near the $k$-core threshold of a random graph}

\maketitle

\begin{abstract}  
 We prove that $G_{n,p=c/n}$ \whp\ has a $k$-regular subgraph if $c$ is at least  {$e^{-\Theta(k)}$ above the threshold for the appearance of a subgraph with minimum degree at least $k$; i.e. an non-empty $k$-core.}  In particular, this pins down the   {threshold for the } appearance of a $k$-regular subgraph to a window  {of  size $e^{-\Theta(k)}$.}
\end{abstract}

In this paper, we study the threshold for  {the Erd\H{o}s-R\'enyi random graph model} $G_{n,p=c/n}$ to have a $k$-regular subgraph where $k$ is fixed. This problem was first studied by Bollob\'as, Kim and Verstra\"{e}te~\cite{bkv} who proved, amongst other things, that  $G_{n,p=c/n}$ w.h.p.\footnote{A property is said to hold {\em with high probability (\whp)} if it holds with probability tending to one as $n\rightarrow\infty$.} has a $k$-regular subgraph when $c$ is at least roughly $4k$.   Letzter~\cite{let} proved that this threshold is sharp\footnote{Meaning that there is a function $\r_k(n)$ 
such that for any $\e>0$, $G_{n,p=c/n}$ w.h.p.\ has no $k$-regular subgraph for $c=\r_k(n)-\e$ and \whp\ has a $k$-regular subgraph for $c=\r_k(n)+\e$.}.

This problem is reminiscent of the $k$-core, a maximal subgraph with minimum degree at least $k$.  Pittel, Spencer and Wormald~\cite{Pittel99} established the  threshold for $G_{n,p}$ to have a {non-empty} $k$-core to be a specific constant $c_k=k+o(k)$ (we specify $c_k$ more precisely in~(\ref{eq:ck}) below).  This provides a lower bound on the threshold for a $k$-regular subgraph, and it is natural to ask:

{\bf Question:} Is the threshold for a $k$-regular subgraph equal to the $k$-core threshold?

Bollob\'as, Kim and Verstra\"{e}te~\cite{bkv} proved that the answer is ``No'' for $k=3$ and conjectured that it is ``No'' for all $k\geq 4$.  On the other hand, Pretti and Weigt~\cite{pw} provided a non-rigorous analysis and claimed that it indicates  the answer is ``Yes'' for $k\geq 4$.

Pra\l{}at,  Verstra\"ete, and Wormald~\cite{Pawel} proved that \whp\ the $(k+2)$-core of $G_{n,p}$ (if it is non-empty) contains a $k$-regular spanning subgraph. Chan and Molloy~\cite{Mike} proved the same for the $(k+1)$-core. So the $k$-regular subgraph threshold is at most $c_{k+1}\approx c_k+1$.  We will reduce this bound to within an exponentially small distance (as a function of $k$) from $c_k$:

\begin{theorem}\label{mt}
For $k$ a sufficiently large constant, and for any $c\geq c_k+e^{-k/300}$, $G_{n,p=c/n}$ \whp\ contains a $k$-regular subgraph.
\end{theorem}

It is not hard to see that the $k$-core cannot have a $k$-regular {\em spanning} subgraph; for example \whp\ it has many vertices of degree $k+1$ whose neighbours all have degree $k$.  Our approach is to start with the $k$-core and repeatedly remove such vertices, along with other problematic vertices.  We will then apply a classic theorem of Tutte
to show that what remains has a spanning $k$-regular subgraph.  The aforementioned papers~\cite{Pawel, Mike} applied Tutte's theorem to the $(k+2)$- and $(k+1)$-core.  

\mmm{
The $k$-core is well known to have size  \mmmm{$(1-o_k(1))n$}, and we will show that we only remove $o_k(1)n$ vertices (see Remark~\ref{rokn}).  So the $k$-regular subgraph that we obtain will have size \mmmm{$(1-o_k(1))n$}.
}

The new arguments required in this work are (i) stripping the $k$-core down to something to which Tutte's theorem can be applied and (ii) applying Tutte's theorem to it.  {The first part requires a delicate variant of the configuration model (see the discussion at the beginning of Section~\ref{stsp}),} whereas the presence of degree $k$ vertices brings new challenges to the second part.

The number of problematic vertices, as described above, is linear in $n$. Furthermore, removing them from the $k$-core will cause a linear number of vertices to have their degrees drop below $k$.  It is not surprising that if  $c$ is too close to $c_k$, then \whp\ what remains will have \mm{an empty} $k$-core, and so this argument will not work unless $c$ is bounded away from $c_k$.   Fortunately, when $k$ is large, the number of problematic vertices is very small (but linear in  $n$):  $e^{-\Theta(k)} n$.  So we only need $c$ to be bounded away from $c_k$ by an exponentially small distance {(in terms of $k$)}.   Furthermore, the subgraph that we show to have a $k$-factor consists of all but $e^{-\Theta(k)}n$ vertices of the $k$-core \mmmm{(see Remark~\ref{rokn}).} {This is consistent with a result of}  Gao~\cite{Gao} {who} proved that any $k$-regular subgraph must contain all but at most $\e_k n$ vertices of the $k$-core where $\e_k\rightarrow 0$ as $k$ grows.

 { {\em Organization of the paper:}  We begin, in section 1, by presenting Tutte's condition.  Section 2 contains some brief probabilistic tools. The stripping procedure to find our subgraph is given and analyzed in section 3; this is most of the work.  Finally, in section 4, we show how to prove that the subgraph satisfies Tutte's condition and thus has a spanning $k$-regular subgraph.
}

\section{Tutte's condition}
We begin by presenting Tutte's theorem for establishing that a graph has a $k$-regular spanning subgraph. Recall that a $k$-regular spanning subgraph is called a {\em $k$-factor}.
Let $\G$ be a graph with minimum degree at least $k$.

\begin{definition}\label{dlh} 
$L=L(\G)$ is  the set of {\em low vertices} of $\G$, i.e.\ the vertices $v$ with $d_{\G}(v)=k$, and $H=H(\G)$ is the set of {\em high vertices} of $G$, i.e.\ the vertices $v$ with $d_{\G}(v)\geq k+1$.   For any set of vertices $Z$, we use $Z_L,Z_H$ to denote $Z\cap L$, respectively $Z\cap H$.
\end{definition}

\mm{{\bf Notation:} } For any $S \subseteq V(\G)$ we use $e(S)$ to denote the number of edges of $\G$ with both endpoints in $S$.
For any disjoint $S,T \subseteq V(\mm{\G})$, we use $e(S,T)$ to denote the number of edges of $G$ from $S$ to $T$ and $q(S,T)$ to denote the number of components $Q$ of $\mm{\G} \setminus (S\cup T)$ such that $k|Q|$ and $e(Q,T)$ have different parity.   \mm{Throughout the paper, we use $d_A(v)$ to denote the number of neighbours that $v$, a vertex, has in $A$, a subset of the vertices.} \dm{Furthermore, we refer to the {\em total degree} of $S$ as the sum of degrees of the vertices in $S$.}

\begin{theorem}[\cite{Tutte}]\label{thm:Tutte} 
A graph $\G$ with minimum degree at least $k\geq1$ has a $k$-factor if and only if for every pair of disjoint sets $S,T \subseteq V(\G)$, 
$$
k|S| \ge q(S,T) + k|T|- \sum_{v \in T} d_{\G \setminus S}(v).
$$
\end{theorem}

\begin{corollary}\label{Tutte-condition}
A graph $\G$ with minimum degree at least $k\geq 1$ has a $k$-factor if \mm{and only if} for every pair of disjoint sets $S,T \subseteq V(\G)$,
\begin{equation}\label{etutte2}
k|S|+\sum_{v\in T_H}(d_{\G}(v)-k) \ge q(S,T)+e(S,T).
\end{equation}
\end{corollary}

{\proofstart}
Rearranging the terms in Theorem~\ref{thm:Tutte}, we obtain that the condition given there for the existence of a $k$-factor is equivalent to  every pair of disjoint sets $S,T \subseteq V(\G)$ satisfying:
\begin{eqnarray*}
\sum_{v \in T} d_{\G}(v) + k|S| &\ge &q(S,T)+k|T|+e(S,T)\\
\mbox{i.e.} \qquad \sum_{v \in T} (d_{\G}(v)-k) + k|S| &\ge &q(S,T)+e(S,T)
\end{eqnarray*}
which is equivalent to (\ref{etutte2}) since $d_{\G}(v)=k$ for all $v\in T_L$.
{\proofend}

In all but one case, we will in fact show that $S,T$ satisfy the stronger condition, which implies~(\ref{etutte2}) since $d_{\G}(v)\geq k+1$ for all $v\in T_H$:
\begin{equation}\label{etutte}
k|S|+|T_H| \ge q(S,T)+e(S,T).
\end{equation}

{\bf Remark:}  In previous papers~\cite{Mike,Pawel} Theorem~\ref{thm:Tutte} was applied to subgraphs with minimum degree at least $k+1$, specifically the $(k+1)$-core~\cite{Mike} and the $(k+2)$-core~\cite{Pawel}. So it sufficed to prove the weaker bound $k|S|+|T| \ge q(S,T)+e(S,T)$.

We begin by showing that in Corollary~\ref{Tutte-condition} we may assume that $ S \subseteq H$ and every component counted by $q(S,T)$ has a high vertex:

\begin{lemma}\label{lem:minimality}
A graph $\G$ with minimum degree $k\geq 1$ has a $k$-factor if \mm{and only if}~(\ref{etutte2}) holds for every pair of disjoint sets $S,T \subseteq V(\G)$
satisfying:
\begin{itemize}
\item [(M1)] $S \subseteq H$; and
\item [(M2)] every component $Q$ counted by $q(S,T)$ satisfies $Q_H\neq\emptyset$.
\end{itemize}
\end{lemma}

\proofstart 
 {We will prove that if~(\ref{etutte2}) holds for every $S,T$ satisfying (M1) and (M2) then~(\ref{etutte2}) holds for every $S,T$.  So Corollary~\ref{Tutte-condition}  implies that $\G$ has a $k$-factor.
To do this, we show that if $S,T$ violate~(\ref{etutte2}) and violate either (M1) or (M2) then we can modify $S,T$ so that (M1) and (M2) both hold but~(\ref{etutte2}) is still violated. \mm{This proves} our lemma.}

Suppose there exist two disjoint {sets} $S,T \subseteq V(\G)$  {that violate both~(\ref{etutte2}) and (M1). So there exists some $u \in S_L$.} We will show that after moving $u$ to $V(\G) \setminus (S \cup T)$,  {(\ref{etutte2})} will still fail for the new pair of sets. By applying this procedure iteratively for every $u \in {S_L}$,  {we obtain two sets violating~(\ref{etutte2}) and satisfying (M1)}.

Let $S':=S \setminus \{u\}$. Note first that $k|S'|=k|S|-k$ and $e(S',T)=e(S,T)-d_T(u)$. Now, since $d_{\G}(u)=k$, there are at most $(k-d_T(u))$ neighbours of $u$ in $V(\G) \setminus (S \cup T)$ (observe that some neighbours of $u$ might be in $S$). In the worst case, all neighbours of $u$ in $V(\G) \setminus (S \cup T)$ belong to different components all contributing to $q(S,T)$; moreover, after moving $u$ to $V(\G) \setminus (S \cup T)$, they form one connected component of $V(\G) \setminus (S \cup T)$, and so do not contribute to $q(S',T)$ anymore. In any case, $q(S',T) \ge q(S,T) - (k - d_T(u))$. Since $|T_H|$ is left unchanged by switching from $S$ to $S'$, we have
\begin{align*}
k|S'|+ {\sum_{v\in T_H}(d_{\G}(v)-k)}&=k|S|+ {\sum_{v\in T_H}(d_{\G}(v)-k)}-k \\
&< q(S,T)+e(S,T)-k \\
& \le q(S',T)+(k-d_T(u)) +e(S',T)+d_T(u) -k \\
&= q(S',T)+e(S',T),
\end{align*}
and hence~(\ref{etutte2}) still fails.

Now suppose  there exist two disjoint $S,T \subseteq V(\G)$  {violating~(\ref{etutte2})} and which satisfy (M1) but  {not} (M2).  Let $Q$ be a component of $\G\bk (S\cup T)$ such that $k|Q|$ and $e(Q,T)$ have different parity and with $Q_H=\emptyset$.  Note that $e(Q, V(\G)\bk Q)=e(Q, S \cup T)$ has the same parity as $\sum_{u\in Q} d(u)=k|Q|$ since  {$Q_H=\emptyset$}.  So $e(Q, V(\G)\bk Q)$ has a different parity than $e(Q,T)$ and thus $e(Q,S)\neq 0$.   Now we move $Q$ into $T$; i.e.\ we set $T'=T\cup Q$.  Since $Q\subseteq L$, we have $T'_H=T_H$.  Since $Q$ is a component of $\G\bk (S\cup T)$, this move does not affect whether any other component counts towards $q$; i.e.\ $q(S,T')= q(S,T)-1$. Moreover, as we argued above, $e(Q,S)>0$ and so $e(S,T')\geq e(S,T)+1$.  Combining this with $k|S|+ {\sum_{v\in T_H}(d_{\G}(v)-k)} < q(S,T)+e(S,T)$ yields $k|S|+ {\sum_{v\in T'_H}(d_{\G}(v)-k)} < q(S,T')+e(S,T')$ and hence~(\ref{etutte2}) still fails. Clearly, $S$ has not changed and hence (M1) still holds. Repeated applications result in a pair $S,T$ that violates~(\ref{etutte2})  and satisfies (M1) and (M2).  
\proofend

\section{Probabilistic preliminaries}\label{spp}

\mmm{We use $\Bin(\ell,p)$  to denote the binomial random variable with $\ell$ trials and success probability $p$.  We use $\Po(x)$ to denote the Poisson variable with mean $x$.}

Pittel, Spencer and Wormald~\cite{Pittel99} established  the $k$-core threshold to be:
\begin{equation}\label{eck}
c_k=\min_{x>0}\frac{x}{1-e^{-x}\sum_{i=0}^{k-2}\frac{x^i}{i!}}.
\end{equation}
In~\cite{Pawel} the asymptotic value of $c_k$ is determined up to an additive ${O(1/\log k) =} o_k(\cdot)$ term.  Setting $q_k=\log k-\log (2\pi)$, we have
\begin{equation}\label{eq:ck}
c_k =k+(kq_k)^{1/2}+\left( \frac {k}{q_k} \right)^{1/2} + \frac{q_k-1}{3}+ O \left( \frac {1}{\log k} \right).
\end{equation}

We will use the following well-known bounds on tail probabilities known as Chernoff's bound (see, for example,~\cite{JLR}, Theorem~2.1). 
\mmm{Let $X$ be distributed as $\Bin(\ell,p)$, so} $\E[X]=\mu=p\ell$.  Then, 
\begin{eqnarray}\label{chernoff:low}
    \Pr{ (X \le \mu - t)} \le \exp \left( -\frac{t^2}{2\mu} \right) 
   \end{eqnarray} 
and
\begin{eqnarray}\label{chernoff:up}
   \Pr{ (X \ge \mu + t) } \le \exp\left(-\frac{t^2}{2(\mu+t/3)}\right). 
   \end{eqnarray}
In addition, all of the above bounds hold for the general case in which $X=\sum_{i=1}^{\ell} X_i$ and $X_i$ is the Bernoulli random variable with parameter $p_i$ with (possibly) different $p_i$'s.

The \mm{Azuma-Hoeffding} inequality can be generalized to include random variables close to martingales. One of our proofs, proof of Lemma~\ref{lstop}, will use the supermartingale method of Pittel et al.~\cite{Pittel99}, as described in~\cite[Corollary~4.1]{Wormald-DE}. 
Let $G_0, G_1,\ldots, G_\ell$ be a random process and let $X_i$ be a random variable determined by $G_0, G_1, \ldots, G_i$, $0 \le i \le \ell$. Suppose that for some \mm{real constant $b$ and real positive constants $c_1,\ldots ,c_{\ell}$}, 
$$
\ex (X_i - X_{i-1} | G_0, G_1, \ldots, G_{i-1}) < b \quad \text{ and } \quad |X_i - X_{i-1}| \le c_i
$$ 
for each $1 \le i \le \ell$. Then, for every $\alpha > 0$,
\begin{equation}\label{eq:HA-inequality2}
\pr \left[ \text{For some } i \text{ with } 0 \le i \le \ell : X_i - X_0 \ge ib+\alpha \right] \le \exp \left(-\frac{\alpha^2}{2 \sum_j c_j^2} \right).
\end{equation}

\mm{Throughout the paper, we often omit floor and ceiling signs.}

\section{Finding the subgraph}\label{sfts}

{Let $k \in \N$} and set
\begin{equation}\label{ebeta}
\b= e^{-k/200}.
\end{equation}
We begin with a random graph $G=G_{n,p}$ with $p=c/n$ for some constant $c$ satisfying 
\begin{equation}\label{eq:conditions_for_c} 
c_k+ k^{10} \beta =: \cmin \le c \le \cmax := c_k + k^{-1/2},
\end{equation}
where $c_k$ is the threshold of the emergence of the $k$-core.
Our goal is to find {(for $k$ sufficiently large)} a subgraph $K$ of the $k$-core with certain properties, which will ensure that it has a $k$-factor.  

 Our first property  is simply a degree requirement. Of course, $K$ must have minimum degree at least $k$; for technical reasons, it will help if the maximum degree is bounded by a constant; we arbitrarily chose {this to be} $2k$.  In the introduction, we noted that $K$ cannot have any vertex of degree greater than $k$ whose neighbours all have degree $k$.  It is not hard to build similar problematic local structures; it turns out that we can eliminate all of them by not allowing any vertex of degree greater than $k$ that has many neighbours of degree $k$;  our second property enforces this.   Our third property  simply says that $K$ has linear size.  Our final property is a trivial necessary condition for having a $k$-factor.  

\begin{itemize}
\item [(K1)] for every vertex $v\in K$,  $k \le d_K(v) \le 2k$;
\item [(K2)] for every vertex $v\in K$ with $d_K(v)\geq k+1$, we have $|\{w \in N_{{K}}(v): d_K(w)=k\}| \le \frac{9}{10}k$;
\item [(K3)]  $|K|\geq \frac{n}{3}$;
\item [(K4)]  $k|K|$ is even.
\end{itemize}
(In fact, we will be able to find an induced subgraph $K$ of $G$ satisfying these properties.)

\subsection{The stripping procedure}\label{stsp}

{In order to achieve our goal, we are going to use a carefully designed stripping procedure during which one vertex is removed in each step  until a subgraph $K$ satisfying properties (K1), (K2), and~(K3) remains. It will be easy to modify $K$ to enforce the final property (K4), if necessary, at the end.} 

 {We found property (K2) to be particularly challenging to enforce.  The typical approach to this sort of problem is to repeatedly remove a vertex if it violates one of (K1-3).  Often one can argue that at every step, the remaining graph  is uniformly random conditional on its degree sequence (for example, this happens when analyzing the $k$-core stripping process).  In some situations, the vertex set is initially partitioned into a fixed number of parts, and one must condition on the number of remaining neighbours each vertex has in each part; this is more complicated but in principle not much more difficult than conditioning on the degree sequence.   In the present situation, enforcing \mm{(K2)} requires conditioning on the number of remaining neighbours each vertex has in $W$, the set of vertices of degree $k$.  However, this is not an initial partition; $W$ changes during the process.  This made our analysis more difficult.

In dealing with this problem, it helps to partition $W$ into $W_0$, the vertices that initially have degree  $k$,  and $W_1$, the vertices whose degrees change to $k$ during the process.  $W_1$ is much smaller than $W_0$ and so we can afford to delete vertices if they have at least {\em two} neighbours in $W_1$ rather than at least $\frac{9}{10}k$.  This simpler deletion rule helps us deal with the fact that $W_1$ is changing throughout our stripping process.}

We begin with the $k$-core of $G$, as any {subgraph} $K$ satisfying (K1) must be a subgraph of the $k$-core. The $k$-core of a graph can be found by repeatedly deleting vertices of degree less than $k$ from the graph so this  {can be viewed as an initial phase of our stripping procedure.} We continue stripping the graph (as explained below) and, throughout our procedure, we partition the vertex set as follows: 
\begin{eqnarray*}
W_0&=&  \mbox{ the vertices in the remaining graph that had degree $k$ in the $k$-core of $G$;}\\
W_1&=& \mbox{ the vertices of degree at most $k$ in the remaining graph that are not in $W_0$;}\\
R&=& \mbox{ the vertices of degree greater than $k$ in the remaining graph.}
\end{eqnarray*}
Note that vertices may move from $R$ to $W_1$ during our procedure, but no vertex leaves $W_0$ \mm{or $W_1$} unless it is deleted.

The following definition governs the stripping procedure.

\begin{definition}\label{def:deletable}
We say a vertex $v$ is {\em deletable} if in the initial $k$-core:
\begin{itemize}
\item [(D1)] $\deg(v)>2k$; \dm{or}
\item [(D2)] $v\notin W_0$ (that is, $\deg(v) \ge k+1$) and $v$ has at least $\hf k$ neighbours in $W_0$;
\end{itemize}
or if in the remaining graph:
\begin{itemize}
\item [(D3)] $\deg(v)<k$; \dm{or}
\item [(D4)] $v\in R$ and $v$ has at least two neighbours that are  in $W_1$; or
\item [(D5)] $v\in W_1$ and $v$ has a neighbour that is either (i) in $R$ and deletable, or (ii) in $W_1$.
\end{itemize}
Furthermore, 
\begin{itemize}
\item [(D6)] once a vertex becomes deletable it remains deletable. 
\end{itemize}
\end{definition}

\noindent {\bf Remarks:} {Let us make three remarks.}

\begin{itemize}
\item [(a)] Deleting vertices in $W_1$ with non-deletable neighbours in $W_1$  is not required for properties (K1), (K2), and~(K3); we only delete them because it helps with our analysis. 

\item [(b)] In many \mm{similar} stripping processes (for example, the $k$-core process), we have the property that  the subgraph we eventually obtain does not depend on the order in which vertices are deleted.  That is not true for our procedure---whether a vertex ever becomes deletable can depend on the deletion order.  However, our goal is only to obtain a subgraph with the desired properties (K1), (K2), and (K3) and so this works for our purpose.

\item [(c)] Deleting a vertex that is deletable may cause some non-deletable vertices to become deletable which, in turn, might force more non-deletable vertices to change their status. However, this ``domino effect'' will eventually stop (possibly with all remaining vertices marked as deletable) because of (D6).
\end{itemize}

At any point of the algorithm, we use $Q$ to denote the set of deletable vertices \mm{that have not yet been removed}. Recall from~(\ref{ebeta}) that $\b= e^{-k/200}$. We will show that \whp\ we reach $Q=\emptyset$ within $\b n$ steps.  It will be convenient to force our stripping procedure to halt after $\b n$ steps regardless of whether it has reached {the desired property (that is, $Q=\emptyset$)}.

\bigskip

Now, we are ready to introduce our stripping procedure. \mm{This procedure can be applied to any graph $G$. Much of the work in this paper will be to analyze what happens when it is applied to $G=G_{n,p}$.}

\noindent{\bf STRIP}
\begin{enumerate}
\item Begin with the $k$-core {of $G$},  and initialize $Q$ to be \mm{the set of} vertices $v$ with $\deg(v)>2k$ or $v\notin W_0$ and $v$ has at least $\hf k$ neighbours in $W_0$; \mm{i.e.\ for which (D1) or (D2) hold.}
\item Until $Q=\emptyset$ or until we have run $\b n$ iterations, let $v$ be the next vertex in $Q$, according to a specific fixed vertex ordering. Let \mmm{$N(v)$ be the set of neighbours of $v$ remaining at this point.}
\begin{enumerate} 
\item Remove $v$ from the graph (and from $Q$). 
\item If any $u \in \mmm{ N(v)\cap R}$ now has degree at most $k$, then move $u$ from $R$ to $W_1$.
\item If any vertex $w\notin Q$ is now deletable, place $w$ into $Q$.
\end{enumerate}
\end{enumerate}

\noindent{\bf Clarification:}  In step 2c, $w$ does not leave whichever of $W_0,W_1,R$ it was in.  So \mm{ for example, $w$ could be in} $Q\cap W_0$.

\noindent{\bf Remark:} 
Note that initially no vertex has degree less than $k$ and $W_1 = \emptyset$ so, indeed, all initially deletable vertices are added to $Q$ in step 1.

\medskip

{The following observation is an immediate consequence of Definition~\ref{def:deletable}. Indeed, parts (a) and (b) follow from (D5); part (c) follows from (D4).}

\begin{observation}\label{ob1}
The following properties hold \mmm{at the beginning of every iteration}.
\begin{enumerate}
\item[(a)] If $u\in W_1$ then $u$ has no neighbours in $W_1\bk Q$.
\item[(b)] If $u\in R\cap Q$  then $u$ has no neighbours in $W_1\bk Q$.
\item[(c)] If $u\in R\bk Q$  then $u$ has at most one neighbour in $W_1$.
\end{enumerate}
\end{observation}

{Here are two more straightforward but useful observations.}

\mmm{
\begin{observation}\label{oq4k2}  During any iteration of Step 2, at most $4k^2$ vertices enter $Q$.
\end{observation}

{\bf Proof:} This follows by noting: (i) the deleted vertex $v$ has at most $2k$ neighbours $u$; (ii) if $u$ enters $W_1$ then $u$ has at most $k-1$ other neighbours $z$; (iii) if a vertex $w$ becomes deletable then it is either one of the $u,z$ mentioned in (ii) or $w\in W_1$ and $w$ is a neighbour of a $z\in N(u)\cap R$ which becomes deleteable; (iv) each such $z$ has at most one neighbour $w\neq u$ that is in $W_1$, else $z$ would have already been in $Q$.
}


\begin{observation}\label{ob2}
If {the stripping procedure} terminates with $Q=\emptyset$, then in the remaining subgraph: the degree {of every vertex} is in $\dm{\{k, k+1, \ldots ,2k\}}$ and {each} vertex in $R$ has {at most} $\mm{\hf k}$ neighbours of degree $k$ {(provided that $k \ge 3$)}.  Thus, the remaining subgraph satisfies {properties} (K1) and (K2).
\end{observation}

\subsection{Configuration models}\label{scm}

We model the $k$-core of $G=G_{n,p=c/n}$ with the configuration model introduced by Bollob\'as~\cite{Bollobas}, and inspired by  Bender and Canfield~\cite{Bender}. We are given the degree sequence of a graph (that is, the degree of each vertex). We take $\deg(v)$ copies of each vertex $v$ and then choose a uniform pairing of those vertex-copies.  Treating each pair as an edge gives a multigraph.  It is well-known (see, for example, the result of McKay~\cite{McKay}) that the multigraph {(for a degree sequence meeting some mild conditions; these conditions are met in our application)} is simple with probability tending to a positive constant, and it follows that if a property holds \whp\ for a uniformly random configuration then \mm{it} holds \whp\ for a uniformly random simple graph on the same degree sequence; see the survey of Wormald~\cite{Wormald_models} for more on this, and for a history of the configuration model and other related models.

Throughout our analysis, we often refer to a pair in the configuration as an edge in the corresponding multigraph.  A sub-configuration is simply a subset of the pairs in a configuration, and thus yields a subgraph of the multigraph.

We define $C$ to be the $k$-core of $G_{n,p=c/n}$. 
We expose the degree sequence $\cald$ of $C$  and then define $\La$ to be a uniform configuration  with degree sequence $\cald$. We will prove:

\begin{lemma}\label{mtc}
W.h.p.\ STRIP terminates with $Q=\emptyset$ when run on $\La$.   
\end{lemma}

As described above, we can transfer our results on random configurations to random simple graphs thus obtaining:

\begin{lemma}\label{mtsg}
W.h.p.\ STRIP terminates with $Q=\emptyset$ when run on $C$.   
\end{lemma}

\subsection{$k$-core properties}

We will require the following properties of the configuration $\La$.

{\bf Setup for Lemma~\ref{lw0}:}  $k$ is a sufficiently large constant, and  $c_k < c \le \cmax = c_k + k^{-1/2}$.  $C$~is the $k$-core of $G_{n,p=c/n}$.  $\La$ is a uniform configuration with the same degree sequence as $C$.   {Finally,} $W_0,R$ are as defined in Section~\ref{stsp}.

\begin{lemma}\label{lw0} \mm{W.h.p.\ } before the first iteration of STRIP:
\begin{enumerate}
\item[(a)] $|\La|> {0.99} n$;
\item[(b)] ${ 0.99} \frac{n}{k} < |W_0| < {1.01} \frac{n}{k}$;
\item[(c)] the total degree of \mm{the set of} vertices with  degree greater than $2k$ is at most {$e^{-k/6}n$};
\item[(d)] there are at least $\frac{n}{{5}k}$ edges with both endpoints in $W_0$;
\item[(e)] there are at least ${\frac {1}{2}} n$ edges from $W_0$ to $R$;
\item[(f)] there are at least $\inv{3}k n$ edges with both endpoints in $R$;
\item[(g)] $C$ has at most {$e^{-k/3}n$} vertices of degree at most $2k$ and with at least $\hf k$ neighbours in $W_0$;
\item[(h)] at least $\frac{n}{{200}}$ vertices in $R$ have no neighbours in $W_0$.
\end{enumerate}
\end{lemma}

\mmmm{The proof of Lemma~\ref{lw0} is straightforward, but lengthy.  So we defer the proof to Section~\ref{sec:lw0}.}

\begin{corollary}\label{ck3} When STRIP terminates, the remaining subgraph has size at least $\frac{n}{3}$; that is, it  satisfies (K3).
\end{corollary}

\proofstart  This follows from Lemma~\ref{lw0}(a), since ${0.99}n-\b n \ge \frac{n}{3}$, {for $k$ sufficiently large}.
\proofend

\subsection{Stripping a configuration}\label{ssac}
\mm{
We will apply STRIP to the random configuration $\La$. In order to analyze this process, We must carefully track  the information that is exposed.   If the partition $W_0,W_1,R$ were fixed throughout our procedure, then we would simply expose the number of remaining neighbours that each vertex has in each part.  But because $W_1$ and $R$ change during the process, our exposure is more delicate.  We found that the best way to deal with this complication includes exposing many of the remaining edges involving $W_1$ and $R$.}

\begin{definition}\label{def:RW}
Suppose the vertices of a configuration are partitioned into $\calr,\calw_0,\calw_1$. 
We define the {\em RW-information} to be:
\begin{itemize}
\item for each vertex $v\in \calw_0$, $\deg_{\calw_1\cup\calr}(v), \deg_{\calw_0}(v)$;
\item for each vertex $v\in \calw_1\cup \calr$, $\deg_{\calr}(v)$, $\deg_{\calw_0}(v)$, $\deg_{\calw_1}(v)$;
\item all vertex-copy pairs that have one vertex-copy in $\calw_1$ and the other in $\calw_1\cup \calr$.
\end{itemize}
\end{definition}

\mm{
{\bf Remark:} A priori, it is not obvious that this is the information we should expose.  For example,  for each $v\in \calw_0$, it may seem more natural to expose  $\deg_{\calr}(v)$ and $\deg_{\calw_1}(v)$ rather than $\deg_{\calw_1\cup\calr}(v)$.  But this precise set of information is what we needed to make our analysis work.

We restate STRIP here, describing how it runs on a configuration and adding details about how we expose the pairs of the configuration. Recall that pairs of vertex-copies correspond to edges in a graph, so when we remove a vertex-copy we do not \pt{necessarily} remove the actual vertex or any other copies of the vertex. 

Recall the procedure STRIP and the definitions of $W_0,W_1,R$ and \pt{being} deletable from Section~\ref{stsp}.}

\medskip

\noindent{\bf STRIP2}
\begin{enumerate}
\item \mm{Begin with $\La$ and set $\calw_0=W_0,\calw_1=\emptyset,\calr=V(\La)\setminus W_0$.  (Note that these sets are $W_0,W_1,R$ respectively at the beginning of our procedure.) Now expose the RW-information.}
\item Initialize $Q$ to be all vertices $v$ with $\deg(v)>2k$ or $v\notin W_0$ and $v$ has at least $\hf k$ neighbours in $W_0$.
\item Until $Q=\emptyset$ or until we have run $\b n$ iterations, let $v$ be the next vertex in $Q$, according to a specific fixed vertex ordering. 
\begin{enumerate} 
\item Expose the partners of every vertex-copy of $v$ (of course, if they are not already exposed). Let \dm{$N(v)$} be the set of neighbours of $v$ \mmm{remaining at this point.}
\item Remove $v$ from $\La$ (and from $Q$), along with all vertex-copies of $v$ and their partners. 
\item If any $u \in \mmm{N(v)\cap R}$ has its degree decreased to at most $k$, then
\begin{enumerate}
\item move $u$ from $R$ to $W_1$;
\item \mm{expose the vertex-copies of $u$ that have  partners in $W_1\cup R$; for each such vertex-copy, expose its partner.}
\end{enumerate}
\item If any vertex $w\notin Q$ is now deletable, place $w$ into $Q$.
\item \mm{Set $\calw_0=W_0,\calw_1=W_1,\calr=R$ and expose the RW-information of the remaining configuration.}
\end{enumerate}
\end{enumerate}

\mm{
{\bf Remark:} The exposure of the RW-information in step \dm{3}(e) is redundant; that information was in fact exposed during previous steps.  We state step \dm{3}(e) in this way to be explicit about the fact that it is exposed.

We define $W_0(i),W_1(i), R(i)$, and $Q(i)$ to be those vertex sets at the end of  iteration $i$ of STRIP2.
We define $\Psi(i)$ to be the RW-information of the remaining configuration with partition $\calr=R(i), \calw_0=W_0(i)$, and $\calw_1=W_1(i)$;
i.e. the RW-information exposed during step 3\dm{(e)}. $W_0(0), W_1(0), R(0), Q(0), \Psi(0)$ are the sets and RW-information from steps 1 and 2.

\begin{observation}\label{orw} We expose enough information to carry out each step of STRIP2.
\end{observation}

\proofstart
$\Psi(0)$ specifies the vertices of $Q$ in step 2. Now consider iteration $i$ of step 3.
The exposed partners of the vertex-copies of $v$ allow us to loop through the vertices $u$ in Step 3(c), and the RW-information $\Psi(i-1)$ tells us the degree of each such $u$.  

To determine which vertices $w$ become deleteable: We know whether $w$ satisfies (D3) from $\Psi(i-1)$ and the number of vertex-copies of $w$ that were removed in previous steps.  If $w\in R$ now satisfies (D4) then $w$ had 0 neighbours \dm{(or 1, respectively)} in $W_1$ at the end of iteration $i-1$ (we know this from $\Psi(i-1)$, and at least 2 neighbours \dm{(or 1, respectively)} of $w$ moved into $W_1$ during step 3(c) (we know this because the neighbours of those vertices in $R$ were exposed in step 3(c)). If $w\in W_1$ now satisfies (D5) then first note that we have exposed all neighbours of $w$ that are in  $W_1\cup R$, either in Step 3(c) if $w$ entered $W_1$ during this iteration, or in $\Psi(i-1)$ otherwise.  Any such neighbour in $R$ that is now deleteable satisfies (D4); we have already described how we know whether it is deleteable.  This is enough to determine whether $w$ satisfies (D5).
\proofend

{\bf Clarification:} It is important to note that in the definition of RW-information,   $\calr,\calw_0,\calw_1$ are not necessarily set to be the sets $W_0,W_1,R$ at some point during STRIP2; they can be any partition of the vertices.  This arises in the proof of Lemma~\ref{luni} below; in particular, it is the reason that we need to prove Claim 2.

}

Given a particular partition into $\calr,\calw_0,\calw_1$, and RW-information, $\Psi$, 
let $\Omega_{\Psi}$ be the set of all configurations on vertex set $\calr\cup \calw_0\cup \calw_1$ with RW-information $\Psi$; that is, all configurations in which each $v$ has $\deg_{\calr}(v),\deg_{\calw_0}(v), \deg_{\calw_1}(v), \deg_{\calw_1\cup\calr}(v)$ equal to the values prescribed in $\Psi$,  and where the set of  pairs with one copy in $\calw_1$ and the other in $\calw_1\cup \calr$ is as prescribed in $\Psi$. 

\mm{The next lemma allows us to analyze STRIP2 by treating the configuration remaining after iteration $t$ as being uniformly selected from $\Omega_{\Psi(t)}$.

\begin{lemma}\label{luni}
For any $t\geq 0$, and any possible set $\Psi(t)$:  every configuration in $\Omega_{\Psi(t)}$ is equally likely to be the subconfiguration remaining after $t$ iterations.
\end{lemma}

}

\proofstart
Let $H$ be \mm{the subconfiguration induced by the vertex-copies that remain after $t$ iterations of STRIP2}.  Consider any $H'\in\Omega_{\Psi(t)}$, and form $\La'$ from $\La$ by replacing the pairs of $H$ with  the pairs of $H'$.

\noindent {\em Claim 1: Each vertex $v$ has the same degree in both $\La$ and $\La'$. 
}

\noindent {\em Proof:}
Indeed, \mm{$\La,\La'$} differ only on the pairs of $H,H'$ and, by construction, $v$ has the same degree in $H'$ as in $H$. The claim holds.
\proofend

Claim 1 says that $\La,\La'$ have the same degree sequence {$\cald$} and so both are equally likely to be chosen as our initial configuration.
The remainder of this proof will establish that if we apply $t$ iterations of STRIP to $\La'$ we will obtain $H'$.  This will imply that $H,H'$ are both equally likely to be what remains after applying $t$ iterations of STRIP2 to the our initial random configuration.  This {will finish the proof of} the lemma.

We define $H(i), H'(i)$ to be the sub-configurations remaining after applying $i$ iterations of STRIP2 to $\La,\La'$, respectively. So $H=H(t)$ and we wish to show that $H'=H'(t)$.  Note that this does not follow \mm{immediately} from the fact that $H'\in\Omega_{\Psi(t)}$. \mm{$\Psi(t)$ specifies, for example, that each vertex $v\in H$ has the same number of neighbours from $\calr=R(t)$ in both $\La$ and $\La'$ at the beginning of the procedure. But it is not obvious that, after running $t$ iterations of STRIP2, $\calr$ is the set of remaining vertices with degree greater than $k$. Claim 2, below, argues that this is indeed the case.}

We let $W_0'(i), W_1'(i), R'(i), Q'(i)$ denote the vertex sets $W_0,W_1,R,Q$ after applying $i$ iterations of STRIP2 to $\La'$.  
In what follows, we use $\deg(v), \deg_R(v)$, etc.\ to denote degrees in $\La$ and $\deg'(v), \deg'_R(v)$, etc.\ to denote degrees  in $\La'$.


\noindent{\em Claim 2: for every $0\leq i\leq t$, $W_0(i)=W_0'(i), W_1(i)=W_1'(i),  R(i)=R'(i)$, and $Q(i)=Q'(i)$.} 

{\em Proof:} By definition, $W_1(0)=W_1'(0)=\emptyset$. Claim 1 implies that $W_0(0)=W_0'(0)$ and $R(0)=R'(0)$.   To prove that $Q(0)=Q'(0)$ we apply Claim 1 and we also need to argue that each vertex $v$ has the same number of neighbours in $W_0(0){=W_0'(0)}$ in both $\La$ and $\La'$.    Note that $W_0(t)\subseteq W_0{(0)}$.  The number of  neighbours $v$ has  in $W_0(t)$ is specified by $\Psi(t)$ to be the same in both $\La,\La'$, and the set of edges from $v$ to $W_0{(0)}\bk W_0(t)$ is identical in both $\La,\La'$ as $W_0{(0)}\bk W_0(t)$ is not in $H$.

Having established the base case, we proceed by induction {on $i$}.  
Suppose that Claim 2 holds for some $i<t$.  Since $Q(i)=Q'(i)$, the \mm{$(i+1)$-st} vertex deleted is the same for both $\La,\La'$, since we choose the next deletable vertex according to the same ordering in both procedures; let $v$ be that vertex.
So the vertex set of the subgraph after $i+1$ steps is the same in both procedures. Furthermore,  $v$ \mm{is not in} $H$ as it is removed during the first $t$ iterations on $\La$. \mm{Thus} the set of pairs including a copy of  $v$ is the same in both $\La,\La'$ and hence in both remaining configurations.   

The key observation is that all decisions made in STRIP2 are based on sets  of pairs that are equal in $H(i), H'(i)$.  This implies that we will make the same changes to the various vertex sets in both configurations.   Indeed, the decisions made are determined entirely by:
\begin{itemize}
\item Pairs containing copies of $v$: these are the same in $H(i),H'(i)$ since $v\notin H$.
\item \mm{The decision as to whether $u$ is moved to $W_1$  in step 3(c) is determined by $\deg(u)$.  Every vertex has the same degree in what remains of both $\La$ and $\La'$}  by Claim~1 and the fact that all pairs removed during the first $t$ steps of STRIP2 are the same in both $\La,\La'$, as they each include at least one vertex not in $H$.
\item The decision as to whether $w$ enters $Q$ in Step 3(d) is determined by \mm{whether $w$ is deleteable.  Whether $w$ satisfies (D3) is determined by $\deg(v)$ which is the same in both processes as described above.  Whether $w$ is now in $W_0,W_1$ or $R$ was decided in step 3(c) and so is the same in both processes.  Whether $w\in R$ satisfies (D4) is determined by the number of pairs containing a copy of $w$ and a vertex-copy in $W_1$; all such pairs are the same in both $\La$ and $\La'$, either because the pair is removed during the first $t$ iterations, or because it is a pair between $W_1(t)$ and $R(t)\cup W_1(t)$ and hence is specified by $\Psi(t)$.  The argument for whether $w\in W_1$ satisfies (D5) is the same.}
\end{itemize}

 So STRIP2 makes the same decisions, and carries out the same steps during iteration $i+1$ on both $\La$ and $\La'$.  This yields the claim for iteration $i+1$ {and so the proof of the claim is finished}.
\proofend

Therefore STRIP{2} removes the same sequence of vertices for the first $t$ iterations on  $\La$ and on $\La'$, and so $H'$ is what remains after $t$ iterations on $C$.   {This finishes the proof of the lemma as explained above.}
\proofend

Lemma~\ref{luni}  allows us to analyze the configuration remaining after any iteration $t$ of STRIP{2} by taking a uniform member of $\Omega_{\Psi(t)}$.  This is fairly simple as the members of $\Omega_{\Psi(t)}$ all decompose into the union of a few configurations which can be analyzed independently using the configuration model.  Given any three disjoint vertex sets $\calw_0,\calw_1,\calr$, a set of RW-information $\Psi$, and a configuration $\La\in\Omega_{\Psi}$ we define:
\begin{itemize}
\item $\La_{\calw_0}\subset \La_{\Psi}$ - the sub-configuration induced by $\calw_0$;
\item $\La_{\calw_1}\subset \La_{\Psi}$ - the sub-configuration induced by $\calw_1$;
\item $\La_{\calr}\subset \La_{\Psi}$ - the sub-configuration induced by $\calr$;
\item $\La_{\calw_0,\calw_1\cup\calr}\subset \La_{\Psi}$ - the bipartite sub-configuration induced by $\calw_0,\calw_1\cup\calr$;
\item $\La_{\calw_1,\calr}\subset \La_{\Psi}$ - the bipartite sub-configuration induced by $\calw_1,\calr$.
\end{itemize}

Note that $\Psi$ specifies the vertex-copies and pairs of $\La_{\calw_1}, \La_{\calw_1,\calr}$; thus it also specifies which vertex-copies in $\calw_1$ are in $\La_{\calw_0,\calw_1\cup\calr}$.  To select a uniform member of $\La_{\Psi}$, we can select the other three configurations independently; that is, 
\begin{enumerate}
\item For each vertex $v\in \calr$, choose a uniform partition of the vertex-copies of $v$ not paired with copies in $\calw_1$ into those that will be paired with copies in $\calw_0,\calr$, according to $\deg_{\calw_0}(v),  \deg_{\calr}(v)$.
\item For each vertex $v\in \calw_0$, choose a uniform partition of the vertex-copies of $v$ into those that will be paired with copies in $\calw_0,\calw_1\cup\calr$, according to $\deg_{\calw_0}(v), \deg_{\calw_1\cup\calr}(v)$.
\item Choose $\La_{\calw_0}$ by taking a uniform matching on the appropriately selected vertex-copies in $\calw_0$.
\item Choose $\La_{\calr}$ by taking a uniform matching on the appropriately selected vertex-copies in $\calr$.
\item Choose $\La_{\calw_0,\calw_1\cup\calr}$ by taking a uniform bipartite matching on the appropriately selected vertex-copies in $\calw_0,\calr$ and the vertex-copies of $\calw_1$ that are (implicitly) specified by $\Psi$ to be paired with $\calw_0$.
\end{enumerate}

To see that this yields a uniform member of $\La_{\Psi}$, note that there is a bijection between the union of these choices and the members of $\La_{\Psi}$; note also that the number of choices for steps 3,4,5 is independent of the partitions chosen in steps 1,2.

In every use of this model, we will use it to analyze the configuration remaining at a particular point during STRIP2. We will set $\calw_0=W_0, \calw_1=W_1, \calr=R$ and so we will use the notation, for example, $\La_{W_1,R}$.

The crucial parameter in our analysis is the branching parameter that comes from exploring $W_0{(0)}$ with a branching process.  It is well-known that  for any $c>c_k$, $W_0{(0)}$ is subcritical and so this parameter is less than \mm{1} (see eg.~\cite{mrbranch,sato}). We need to show how far away it is from \mm{1} when $c={\cmin}$.  This bound is the only reason that we require $c\geq {\cmin}$ rather than $c>c_k$ in~(\ref{eq:conditions_for_c}).

\begin{lemma}\label{lbr}  
Let $k \in \N$ be a sufficiently large constant and set $\a=k^{9}\b$.  Let $G=G_{n,p}$ be a random graph with $p=c/n$ for some constant $c_k + k\a = \cmin \le c \le \cmax = c_k + k^{-1/2}$. Then, w.h.p.\ in the $k$-core of $G$, i.e.\ at step $i=0$ of STRIP:
\begin{equation}\label{elbr}
 \frac{\sum_{u\in W_0}\deg_{W_0}(u)(\deg_{W_0}(u)-1)}{\sum_{u\in W_0}\deg_{W_0}(u)} < 1-\a.
\end{equation}
\end{lemma}

\proofstart Recall from~(\ref{ebeta}) that $\b=e^{-k/200}$ and so $\a$ is very small.
\mm{We are going to use the notation and observations used in the proof of Lemma~\ref{lw0}.
See (\ref{eq:c}) for the definition of $f(x)$ and $x$, (\ref{eq:xk}) for the definition of $x_k$, and (\ref{eq:deg_distribution}) for the degree distribution of the $k$-core. In particular, recall that $c=f(x)$ and $c_k=f(x_k)$.  Our first step is to
use the fact that $c$ is bounded away from $c_k$, to bound $x$  away from $x_k$. 
We set
\[x=x_k+\d.\]

\noindent{\em Claim:} $\delta >k \alpha$.

\noindent
{\em Proof:} Recall from~(\ref{edbound}) that $c\leq\cmax$ implies $\d\leq\log k$. Also recall that  $\e=o_k(1)$ was defined in~(\ref{eq:error2}) and used in~(\ref{eq:bound_for_f'}). Since $\a<(\log k)/k$, and since $\e+k\a=o_k(1)$, it follows from~(\ref{eq:bound_for_f'}) that over the range $[x_k,x_k+k\a]$, $f'$ does not exceed $(\e+k\a)\sqrt{(\log k)/k}+O(1/\sqrt{k})<\sqrt{(\log k)/k}$ for $k$ sufficiently large.  Therefore, $f(x_k+k\a)< f(x_k) + (k\a)\cdot \sqrt{(\log k)/k}<c_k+k\a=\cmin \leq c$.  Since $c=f(x)=f(x_k+\d)$ and $f$ is increasing, this implies the claim.
}

We consider the configuration $\La$ with the same degree sequence as the $k$-core of $G$. Recall that $W_0$ is the set of vertices with degree exactly $k$ in $\La$. We could determine the degree sequence of the subconfiguration induced by $W_0$ and then bound the LHS of~(\ref{elbr}), but we obtain a simpler calculation by considering  the following experiment:  choose a uniformly random vertex-copy \mm{$\sigma$} from $W_0$ conditional on \mm{$\sigma$} being paired in $\La$ with another copy from $W_0$; let $u\in W_0$ be the vertex of which \mm{$\sigma$} is a copy.  Set $Z$ to be the number of other copies of $u$ that are paired with vertex-copies in $W_0$.  Note that $\Exp(Z)$ is the LHS of~(\ref{elbr}).

One way to choose \mm{$\sigma$} is to repeatedly take a uniform vertex-copy from $W_0$ and expose its partner; halt the first time that the partner is also in $W_0$.  By Lemma~\ref{lw0}(b), \dm{ \whp\ } a linear proportion of the copies are in $W_0$ and so \whp\ we only expose $o(n)$ \mmm{pairs of vertex-copies} before halting.  Now, having found $u$, we expose
the partners of the remaining $k-1$ copies of $u$.  So $\Exp(Z)$ is simply the expected number of these partners that are in $W_0$. \dm{Since a vertex in $W_0$ has degree $k$, $Z \le k$ always holds, and on the other hand $\Exp(Z)=\Omega(1)$, as by Lemma~\ref{lw0}(b), \whp\  $|W_0| \ge 0.99n/k$, and in such case $\Exp(Z)=\Omega(1)$. Hence, the first order term of $\Exp(Z)$ stems from the event when Lemma~\ref{lw0}(b) holds and only $o(n)$} \mmm{pairs of vertex-copies} were exposed before halting. \mmm{In that case, the probability that a particular copy of $u$ selects a partner in $W_0$ is simply the total degree of $W_0$ divided by the total degree of $\Lambda$ plus a $o(1)$ term for the exposed pairs.  Applying~(\ref{eq:deg_distribution}) to obtain this ratio, we find that $E(Z)=(1+o(1))g(x)$ where}

\[
g(x) := (k-1)\cdot  \frac{k \ e^{-x}\ \frac{x^k}{k!}}{\sum_{i \ge k} i\ e^{-x}\ \frac{x^i}{i!}} =\frac{k(k-1)\ e^{-x} \ \frac{x^{k}}{k!}}{x \cdot e^{-x} \sum_{i \ge k-1}\frac{x^i}{i!}} 
= \frac{k(k-1)\ e^{-x} \ \frac{x^{k}}{k!}}{x \cdot  \Pr{ (\Po(x) \ge k-1). }}
\]

As mentioned in the proof of Lemma~\ref{lw0}\mm{(b)}, $x_k$ minimizes the function $f(x)$. \dm{Using the fact that} $f'(x_k)=0$, it is a simple exercise to show that $g(x_k) = 1$; see~\cite{Molloy_cores} for the details. Hence, using the fact that  $\Pr{ (\Po(x) \ge k-1) } \ge \Pr{ (\Po(x_k) \ge k-1) }$ for $x>x_k$, it follows from~(\ref{eq:error1}) that
\begin{align*}
g(x) & \le \frac{k (k-1)\ e^{-x_k}\ \frac{x_k^k}{k!} \left( 1 - \delta (\log k/k)^{1/2} (1+o_k(1)) \right) }{(x_k + \delta) \ \Pr{ (\Po(x_k) \ge k-1) }}\\ 
& = g(x_k) \left( 1 - \delta (\log k/k)^{1/2} (1+o_k(1)) \right) \\
& \le 1 - \frac {1}{2} \delta (\log k/k)^{1/2}< 1-\alpha, \qquad\qquad\qquad\qquad\mbox{ by \mm{our Claim}}
\end{align*}
provided $k$ is large enough.
\proofend

\subsection{The procedure terminates quickly}

In this section, we will prove that \whp\ STRIP2 terminates with $Q=\emptyset$ when run on $\La$.
In order to show that $Q$ reaches $\emptyset$ within $\b n$ iterations, we will keep track of a weighted sum of the total degree of the vertices in $Q$, and we will show that this parameter drifts towards zero.

Within this parameter, the change in the number of edges from $Q$ to $W_0$ is the most delicate. In particular, the most sensitive part of our process is avoiding cascades that could be formed when the deletion of vertices in $W_0\cap Q$ causes too many other vertices in $W_0$ to be added to $Q$.  (Roughly speaking, Lemma~\ref{lbr} ensures that such cascades do not occur.)
So we place a high weight on the number of edges from $Q\cap W_0$ to $W_0$. We place an even higher weight on the  edges from $Q\bk W_0$ to $W_0$; these edges play a different role in the analysis because they initiate the potential cascades.  \mm{We express this weighted sum with the following variables; each refers to the sets $W_0,W_1,R$ at the end of iteration $i$ of STRIP2.}
\begin{eqnarray*}
A_i&=& \sum_{v\in Q\cap W_0} \deg_{W_0}(v)\\
B_i&=& \sum_{v\in Q\bk W_0} \deg_{W_0}(v)\\
D_i&=& \sum_{v\in Q} \deg_{W_1\cup R}(v)\\
X_i&=&A_i +kB_i + \dm{k^7}\b D_i 
\end{eqnarray*}
Note that, {indeed}, since $\b= e^{-k/200} < k^{-10}$, the edges counted by $A_i$ and $B_i$ have much higher weights in $X_i$ than those counted by $D_i$.

We will prove that $X_i$ has a negative drift, and that \whp\ it \mm{reaches} zero before $\b n$ iterations of STRIP2.  At that point, if any vertices remain in $Q$ then they all have degree zero and so will be removed without any new vertices being added to $Q$.

\begin{observation}\label{obxx} \mm{W.h.p.\ } throughout STRIP2, we always have:
\begin{enumerate}
\item[(a)] $|W_1|\leq 3k\b n$;
\item[(b)] $\dm{e(R, W_0)}\geq \frac{n}{3}$;
\item[(c)] $\dm{e(R,R)}\geq \frac{kn}{4}$;
\item[(d)]  $|Q|\leq \mmmm{5k^2}\b n$. 
\end{enumerate}

\end{observation}

\proofstart 
\mmm{\dmm{In a nutshell, t}he observation follows by noting that, since  STRIP2 carries out a very small number of iterations and hence deletes a very small number of vertices,  these parameters will not change much from their initial values as provided  in Lemma~\ref{lw0}.
\dmm{In more detail}, STRIP2 runs for} at most $\b n$ iterations, where $\b=e^{-k/200}$ from~(\ref{ebeta}), and in each iteration we delete one vertex. 

\mmm{Initially, $W_1=\emptyset$ and every vertex that moves to $W_1$ is the neighbour of a deleted vertex. So $|W_1|$ is bounded by the total degree of the deleted vertices.}
The total degree of the vertices from the initial set $Q$ is at most {$e^{-k/6}n + (2k)e^{-k/3}n < k\b n$} by Lemma~\ref{lw0}(c,g).  Every other deleted vertex has degree at most $2k$ and so the total degree of all deleted vertices is at most $(2k) \b n+k\b n=3k\b n$;
 \mmm{this proves part~(a).} 

\dm{For part~(b), observe that by Lemma~\ref{lw0}(e), initially, \whp\ $e(R,W_0) \ge \frac12 n$. Since by the proof of part~(a) the total degree of all deleted vertices is at most $3k \b n$, and the number of edges deleted is bounded by the total degree of all deleted vertices, at any time $e(R, W_0) \ge \frac12n - 3k\b n \ge \frac13 n$. 
By an analogous argument, this time using Lemma~\ref{lw0}(f), $e(R,R) \ge \frac{nk}{3} - 3k\b n \ge \frac{nk}{4}$, and part~(c) follows.  }

Finally, for part~(d): \mmmm{
Lemma~\ref{lw0}(c,g) implies that the size of the initial set $Q$ is at most $(e^{-k/6}+ e^{-k/3})n<\b n$. Observation~\ref{oq4k2} says that at most $4k^2$ vertices are added to $Q$ during each of the at most $\b n$ steps.  So $|Q|$ can never exceed $\b n+4k^2\b n<5k^2\b n$.} 
\proofend



 The key parameter in bounding the drift of $X$ is {controlled by the following lemmas.}
 
 \mm{If $\s$ is a vertex-copy, then we use $u(\s)$ to denote the vertex that $\s$ is a copy of.}

\begin{lemma}\label{lq0}
\mm{W.h.p.\  at every} iteration $i\leq \b n$ of STRIP2, \mm{the RW-information $\Psi_i$ is such that:}  if a vertex-copy $\s$ is chosen uniformly at random from the remaining $W_0$-copies in $\mm{\La}$, then
\[\ex\left(\deg_{W_0}(u(\s))\right)<\left(2-\frac{\a}{2}\right),\]
where $\a=k^{9}\b$ is from Lemma~\ref{lbr}. 
\end{lemma}

\proofstart  
The expectation we are bounding is {equal to}:
\begin{equation}\label{ed1}
\sum_{u\in W_0}\deg_{W_0}(u)^2/\sum_{u\in W_0}\deg_{W_0}(u).
\end{equation}
By Lemma~\ref{lbr}, at iteration $i=0$ this is 
\[\frac{\sum_{u\in W_0}\deg_{W_0}(u)(\deg_{W_0}(u)-1)+\sum_{u\in W_0}\deg_{W_0}(u)}{\sum_{u\in W_0}\deg_{W_0}(u)} < 1-\a +1 = 2-\a.\]
At iteration $1\leq i\leq \b n$, each of the $i$ vertices that have been deleted decreases the denominator of~(\ref{ed1}) by at most $2k$ (the largest effect is when we remove a vertex of $W_0$ with $k$ neighbours in $W_0$) and so the denominator has decreased by at most $2k\b n$.  The numerator has not increased. Since the denominator was initially at least $n/(5k)$ by Lemma~\ref{lw0}(d), the value of~(\ref{ed1}) is at most
\[
(2-\a) \cdot \frac {n/(5k)}{n/(5k) - 2k\b n} = \frac{2-\a}{1-\frac{2k\b n}{n/5k}}  < (2-\a)(1+10 k^2\b)<2-\frac{\a}{2},
\]
as $\a=k^{9}\b$.
\proofend

This brings up to our key lemma:

\begin{lemma}\label{ldriftx}
\mm{W.h.p.\  at every iteration $i\leq \b n$ of STRIP2, the RW-information $\Psi_i$ is such that:} if the vertex $v\in Q$ that is deleted has degree at least one, then
\[\ex(X_{i}-X_{i-1})\leq -\hf \dm{k^7}\b.\]
\end{lemma}

{\bf Remark:}  If $v$ has degree zero then  $X_{{i}}=X_{{i-1}}$ {(deterministically)} as during that iteration of STRIP2,  no vertex-copies will be removed from the configuration, and no vertices will join $Q$.

\dm{
\proofstart
First note that $X_{i}$ can only increase through vertices joining $Q$.  So most of this analysis focuses on the expected number of vertices that are added to $Q$.  \mm{This analysis relies on the fact that all RW-information is specified deterministically by $\Psi$ (see Definition~\ref{def:RW}), and that the configuration is uniform amongst all configurations with that RW-information (see Lemma~\ref{luni}).}
 
$\Psi$ specifies that vertex $v$ {is incident to} \mm{$\deg_{W_0}(v)$} edges to $W_0$, \mm{$\deg_{W_1}(v)$} edges to $W_{1}$, \mm{$\deg_{R}(v)$} edges to $R$ \mmmm{and/or $\deg_{W_1\cup R}(v)$ edges to $W_1\cup R$ depending on whether $v\in W_1 \cup R$ or $v\in W_0$.} We consider the effect on $X_{{i}}-X_{{i-1}}$ of deleting each of these pairs; \mmm{specifically, the effect of deleting a copy of $v$ and a copy of some $u\in N(v)$ which is specified to be in} \mmmm{$W_0,W_1, R$ or $W_1\cup R$.} 

{\em Case 1: $v\in W_1$.} \dm{Subcase $u \in W_1$:}  The removed pair is not random; it is specified \mm{by $\Psi$, as $\Psi$ specifies all pairs with one member in $W_1$ and the other in $W_1 \cup R$. } By Observation~\ref{ob1}(a), \dm{$u$} is already in $Q$. So no new vertices are added to $Q$ and $D_{i}$ decreases by \mmm{exactly} two. \mmt{Thus,  the deletion of $uv$ causes $X$ to change by}
\[\partial X=-2\dm{k^7}\b.\]

\dm{Subcase $u \in R$}:   The removed pair is not random; it is specified by $\Psi$.  By Observation~\ref{ob1}(c),  \dm{$u$ is} already in $Q$ or has no other neighbours in $W_1$; either way, any neighbours that $u$ has in $W_1$ are already in $Q$. Also, those neighbours are specified by $\Psi$ (since they are edges from $R$ to $W_1$) and so they do not need to be exposed.   

Thus, any new vertices that are added to $Q$ are the result of the neighbours  of $u$ that are in $R$, and so we turn our attention to those neighbours. {Note that the degree of $u$ is at least $k+1$ before the removal of the pair of copies $uv$ (since $u \in R$); if $u$ has degree at least $k+2$ before the removal, no more vertices will be moved to $Q$.}  So for the remainder of the analysis, we will assume that  $u$ has degree $k+1$ \dm{before the removal of $uv$}. Thus $u$ moves to $W_1$ and we expose its at most $k$ neighbours in $W_1\cup R$.   There are two ways that the choice of one such partner, $w$, can cause vertices to be added to $Q$:

\begin{enumerate}
\item[(a)]  \mmt{$w\in R$ and} $w \in Q$.  Then $u$ is added to $Q$ if it was not already in $Q$. \mmt{(Note that if $w\in W_1$ then $u$ is already in $Q$ as it had two neighbours in $W_1$.)}
\item[(b)]  \mmt{ $w\in R$}, $w\notin Q$ and has exactly one neighbour $w'\neq u$ in $W_1$.  So $w$ is added to $Q$, and also $w'$ is added to $Q$ if it was not already in $Q$.  If $w$ has more than one such neighbour $w' \in W_1$, $w'\neq u$, then $w$ and all such neighbours would already be in $Q$. 
\end{enumerate}

\mmt{Both of those situations require $w\in R$, so to analyze the effect of those possibilities we consider the at most $k$ neighbours $u$ has in $R$.}
To expose a neighbour \dm{$w$} of $u$ in $R$, we choose the partner of a \mm{ copy  of $u$} in $\La_R$; that is, we choose a uniform vertex-copy in $\La_R$; $w$ is the vertex containing that copy.

Note that, since $w\in R$, the edges from $w$ to $W_1$ are specified by $\Psi$, and so we do not need to expose any new partners of copies of $w$ to determine whether $w$ has a neighbour $w'\in W_1$. \mmmm{Furthermore, the vertices $w$ that would result in additions to $Q$ as in (a,b) above are specified by $\Psi$.  So to bound the expected change in $X$, we bound the number of copies of such $w$.}
By Observation~\ref{obxx}(a,c,d):\\

\begin{enumerate}
\item[(i)] There are at least $kn/2$ vertex-copies in $\La_R$ to choose from.
\item[(ii)] At most $2k\cdot \mmmm{5k^2} \b n + \b n < \mmmm{11k^3} \b n$ of them are copies of vertices in $R\cap Q$ (as the total degree of  vertices with degree greater than $2k$ is at most ${e^{-k/6} n <} \b n$ by Lemma~\ref{lw0}(c)).  \mmmm{If we choose one of these copies then we have situation (a) and so up to one vertex could be added to $Q$.}
\item[(iii)] At most $6 k^3 \b n$ of them are copies of the at most $3 k^2 \b n$ neighbours of $W_1$ not already in $Q$ (as each such neighbour has degree at most $2k$).  \mmmm{If  we choose one of these copies then we have situation (b) and so up to two vertices could be added to $Q$.}
\end{enumerate}
\mmmm{Recall that in this portion of the analysis, $\deg(u)=k+1$ and  one of $u$'s neighbours, $v$,  has been deleted. } So the expected number of  vertices added to $Q$ is at most $k \cdot \mmmm{23k^3} \b n/(kn/{2}) \leq \mmmm{46k^3}\b$. Each vertex added to $Q$ will increase $X$ by at most $2k^2$ (the extreme case is if it has $2k$ neighbours in $W_0$).
In addition, the removal of the pair $uv$  causes $D_{i}$ to decrease by at least one.  So the expected change to $X$ by the removal of the pair $uv$ is at most:
\begin{equation}\label{eqx1}
\partial X \leq -\dm{k^7}\b + 2k^2\cdot \mmmm{46 k^3}\b <-\hf \dm{k^7}\b,
\end{equation}
{provided $k$ is large enough.}

\dm{Subcase $u \in W_0$}:  We remove a copy of $v$ and its partner \mmm{in $\La_{W_0,W_1\cup R}$; that copy} is selected by choosing a uniform vertex-copy from the $W_0$-copies in $\La_{W_0,W_1\cup R}$.  \mm{This specifies $u$, the vertex to which the copy belongs, and $u$ is added} to $Q$. No other vertices are added to $Q$ in this step. $A_{i}+D_{i}$ can increase by at most $k-1$, as the vertex of $W_0$ added to $Q$ has remaining degree at most $k-1$.  On the other hand, the removal of the pair causes $B_{i}$ to decrease by one.  So \mmmm{ the deletion of $uv$ causes $X$ to change by 
\[\partial X \leq -k + (k-1) = -1. \]
}

{\em Case 2: $v\in R$.} \dm{Subcase $u \in W_1$}:  The removed pair is not random; it is specified by $\Psi$.  By Observation~\ref{ob1}(b), the other endpoint of the edge is already in $Q$. So no new vertices are added to $Q$ and $D_{i}$ decreases by exactly two. \mmmm{ Thus, as in Case 1, the deletion of $uv$ causes $X$ to change by 
\[\partial X=-2\dm{k^7}\b.\]
}

\dm{Subcase $u \in R$:}   \mmmm{There are two differences between this subcase, and the corresponding subcase in Case 1.  
\begin{enumerate}
\item[(1)] The number of vertices $u\in N(v)$ that are in $R$ is specified by $\Psi$, but the actual vertices are not specified by $\Psi$ - they are selected randomly.
\item[(2)] It is possible that $u\notin Q$ and $u$ has a neighbour in $W_1$.  
\end{enumerate} 
}

As in \mmmm{Case 1, subcase $u\in R$,} the degree of $u$ prior to the removal of $uv$ is at least $k+1$. If that degree  at least $k+2$, then no vertex is moved to $Q$ as a result of the deletion of $uv$.  So we assume that   degree  is $k+1$. Thus $u$ moves to $W_1$ and we expose its at most $k$ neighbours in $W_1\cup R$.   There are three ways that the choice of one such partner, $w$, can cause vertices to be added to $Q$.  The first two, arising when $w\in R$, are the same as (a,b) from Case 1, and the analysis of the effect of those possibilites on the expected change in $X$ is the same as in Case 1.  The third is:

\noindent (c) $w\in W_1$.  In that case, $u$ moves to $Q$ and so does $w$ if it is not already in $Q$.  Note that if $u$ has more than one neighbour in $W_1$ then $u$ and those neighbours were already in $Q$.

We expose each of the $\deg_R(v)$ neighbours of $v$ in $R$ by choosing a vertex-copy uniformly from those in $\La_R$; $u$ is the vertex containing that copy.

The set of vertices $u$ whose choice would result in situation (c) is specified by $\Psi$, and there are at most  $|N(W_1)|\leq |W_1|\times k$ such vertices. Each such vertex has at most $k$ copies in $\La_R$, as it has degree $k+1$ and one neighbour in $W_1$, so by Observation~\ref{obxx}(a). there are at most $3k^3\b n$ copies of such vertices in $\La_R$. By Observation~\ref{obxx}(c), there are at least $kn/2$ vertex-copies in $\La_R$ to choose from. So the probability that our choice of $u$ results in situation (c) is at most  \mmmm{$3k^3\b n/(\hf kn)=6k^2\b$}.  If we choose such a copy then up to two vertices will move to $Q$, and each can increase $X$ by  at most $2k^2$.  So the expected impact on $X$ \mmmm{by exposing the vertex $u$} is at most \mmmm{$2k^2 \cdot 6k^2\b = 12k^4\b$.}

\mmt{After exposing $u$, we next expose $N_R(u)$ in order to determine whether any vertices were added to $Q$ because of situations $a,b$.  The same analysis as for~(\ref{eqx1}) shows that the expected impact on $X$ of this step is at most $-\dm{k^7}\b + 2k^2\times {46 k^3}\b$.  As in the calculation for~(\ref{eqx1}), the removal of $uv$ causes $D_i$ to decrease by 1. } Putting this together, the expected change in $X$ resulting from the exposure and removal of the edge $uv$ is at most
\begin{equation}\label{eqx2}
\partial X \leq \mmmm{ -\dm{k^7}\b + 2k^2\cdot {46 k^3}\b + 12k^4\b<-\inv{2}k^7\b},
\end{equation}
provided $k$ is large enough.

\dm{Subcase $u \in W_0$:} The same analysis as in Case 1, \mmm{subcase $u\in W_0$,} shows that regardless of which vertex is selected, the change on $X$ will be $\partial X \leq -1$.

{\em Case 3: $v\in W_0$.}  \mmmm{In this case, $\Psi$  specifies the number of edges from $v$ to $W_1\cup R$, but does not specify how many go to $W_1$ and how many go to $R$.}

\dm{Subcase $u \in W_1\cup R$}: The edge is a pair from 
$\La_{W_0,W_1\cup R}$ and so we choose \mmm{ a vertex-copy uniformly from the $W_1\cup R$-copies in $\La_{W_0,W_1\cup R}$; $u$ is the vertex containing that copy.}  

\mmmm{If we choose $u\in R$ then situations (a,b,c) above are the ways in which this can cause $Q$ to increase.  We can apply the same analysis as for~(\ref{eqx2}).}  The only difference is that the number of $R$-copies in $\La_{W_0,W_1\cup R}$ is at least 
 $n/{3}$ rather than at least $kn/{2}$, as we apply Observation~\ref{obxx}(b) rather than Observation~\ref{obxx}(c).  The result is that the expected change to $X$ conditioning on the selected $u$ being in $R$ is at most 
 \[\partial X \leq \mmmm{ -{k^7}\b + 2k^2\cdot {69 k^4}\b + 18k^5 \b}.\]
 
 If we choose $u\in W_1$, then $u$ will enter $Q$, if it is not already in $Q$; no other vertices will enter $Q$.  
By Observation~\ref{obxx}(a,b) there are at least $n/{3}$ $R$-copies to choose from and at most $3k^2\b n$ $W_1$-copies. So the probability that we select {$u \in W_1$} is at most $3k^2\b n/ (n/{3})={9} k^2\b$. If $u$ is added to $Q$ then at most this will increase $B_i$ by $k-1$ and thus increase $X$ by $k(k-1)<k^2$.
Putting this together, the expected change in $X$ resulting from the removal of the pair $uv$ is at most:
\mmmm{
\[\partial X \leq 9k^2\b\cdot k^2+(1-9k^2\b)\cdot\left( -{k^7}\b + 2k^2\cdot {69 k^4}\b + 18k^5 \b\right)
\leq -\hf k^7\b,\]
} provided $k$ is large enough.

\dm{Subcase $u \in W_0$}:  The edge is a pair from $\La_{W_0}$ and so \mmm{we choose a uniform vertex-copy from  $\La_{W_0}$; $u$ is the vertex containing that copy.}  By Lemma~\ref{lq0}, the expected increase in $A_{i}$ \mmmm{from adding $u$ to $Q\cap W_0$} is at most $2- \frac {\a}{2}$ where $\a=k^{9}\b$. The increase in $D_i$ is at most $k$, and deleting the pair $uv$ decreases $A_{i}$ by two. Putting this together, 
\[\ex(X_{i}-X_{i-1}) \leq -2 + \left( 2- \frac {\a}{2} \right) + k\cdot \dm{k^7}\b < - \frac {\a}{4} < -\hf \dm{k^7}\b.\]

So in every case, the deletion of a copy of $v$ and its partner \dm{$u$} results in an expected change in $X$ of less than $-\hf \dm{k^7}\b$.  \mmmm{Therefore
\[\ex(X_{i}-X_{i-1})\leq -\hf \dm{k^7}\b\cdot \deg(v).\]}
Since $v$ has degree at least one,  this yields the lemma.
\proofend
}

Our bounds on the drift of $X_{i}$ and the initial size of $Q$ imply that our procedure stops quickly.

\begin{lemma}\label{lstop}
W.h.p.\ STRIP2 halts with $Q=\emptyset$ within $\b n$ iterations.
\end{lemma}

\proofstart  We begin by showing that \whp\ we reach $X_i=0$ long before step $i=\b n$.

As we argued in the proof of \dm{Observation}~\ref{obxx}, Lemma~\ref{lw0}(c,g) \mm{implies} that at iteration $i=0$, the total degree of the vertices in $Q$ is \whp\ at most $(e^{-k/6}+2k\cdot e^{-k/3})n$ and so we have $X_0\leq k (e^{-k/6}+2k\cdot e^{-k/3})n<e^{-k/10}n$. Recall from~(\ref{ebeta}) that $\b=e^{-k/200}$.

By Lemma~\ref{ldriftx}, for every $1 \le i \le \b n$, 
\[\E{ (X_{i+1}-X_{i})} \le -b, \qquad \mbox{ where } b:=\frac12 \dm{k^7} \beta.\]
 Note also that (deterministically), for every $1 \le i \le \b n$, $|X_{i+1}-X_i| \le (2k)(4k^2) = 8k^3$, since the degree of a vertex not belonging to $Q$ is at most $2k$, and \mmm{by Observation~\ref{oq4k2} }we add at most $4k^2$ vertices to $Q$ in each iteration.
 
We will use the Martingale inequality from the end of Section~\ref{spp}.  We cannot apply this directly to $X_i$ since $X_i$ stops changing when it reaches zero.  So instead, we couple $X_i$ to a process which is allowed to drop below zero. We define $X_{i+1}'=X_{i+1}$ whenever $X_{i}>0$, and if $X_i=0$ then $X_{i+1}'=X'_i-1$ with probability $b$ and $X'_{i+1}=X_i'$ otherwise.  So for all $i$ we have $\E{ (X'_{i+1}-X'_{i})} \le -b$  and $|X'_{i+1}-X_i'| \le  8k^3$.
Setting $i^* := \mmmm{\frac{\b}{k^3 }}n$, we have
\[\ex(X'_{i^*})\leq  X'_0-i^*b< e^{-k/10}n- \frac12 \mmmm{k^4} \beta^2 n<- e^{-k/50}n.\]
Applying~(\ref{eq:HA-inequality2}) with $\dmm{\alpha=e^{-k/50}n}$, \dmm{$\ell=i^*$ and $c_1,\ldots, c_\ell=8k^3$} yields
\[
\pr(X_{{i^*}}>0) \le \pr(X'_{{i^*}}>0) \leq \exp \left( - \frac {(e^{-k/50} n)^2}{\frac{2\b}{\mmmm{k^3}} n(8k^3)^2} \right) = \exp( - \Omega(n) ) = o(1).
\]
At this point, $|Q| \le e^{-k/10}n+\mmmm{4k^2}i^* <\mmmm{\frac{5\b}{k}}n $, since \mmm{by Observation~\ref{oq4k2}} we add at most $4k^2$ vertices to $Q$ in each iteration. If $X_{i^*}=0$ then there are no edges from $Q$ to the remaining vertices outside of $Q$.   From that point on, no vertices will be added to $Q$ and so STRIP2 halts after $|Q|$ further steps.
So the total number of steps is \whp\ at most $i^*+\mmmm{\frac{5\b}{k}}n<\b n$.
\proofend

This proves Lemma~\ref{mtc} since carrying out STRIP2  performs the same steps as carrying out STRIP on $\La$ (the only difference is that STRIP2 also exposes some information).  Since properties that hold \whp\ on $\La$ also hold \whp\ on $C$, the $k$-core of $G_{n,p=c/n}$, (see the discussion in Section~\ref{scm}), this yields Lemma~\ref{mtsg}.

\subsection{Enforcing (K4)}\label{sk4}
Observation~\ref{ob2}, Lemma~\ref{mtsg}, Corollary~\ref{ck3} and Lemma~\ref{lstop} imply that \whp\ STRIP produces a subgraph $K$ satisfying properties (K1), (K2) and (K3).  It only remains to enforce:

\noindent (K4) $k|K|$ is even.

To do so, we prove that \whp\ the output of STRIP has a vertex $v$ which can be deleted without violating (K1), (K2) or (K3).  Thus, if (K4) does not hold, then we remove $v$ to obtain our desired subgraph $K$.

\begin{lemma}\label{lk4}  W.h.p.\ the output of STRIP contains a vertex of degree greater than $k$ whose neighbours all have degree greater than $k$.
\end{lemma}

\proofstart  We argue that the lemma holds when running STRIP2 on $\La$.  It immediately follows for running STRIP on $C$.  

By Lemma~\ref{lw0}(h),  \whp\ initially at least $\frac{n}{{200}}$ vertices of $R$ have no neighbours in $W_0$.   Since we remove at most one vertex per iteration, at most $\b n$ of these vertices are not in the output.  By Observation~\ref{obxx}(a), the output satisfies $|W_1| \le 3k\b n$.  Each vertex in $W_1$ has degree at most $k$ {(in fact, at this point exactly equal to $k$)}, and so at most $3k^2\b n$ members of $R$ have a neighbour in $W_1$.  This implies the lemma as $\frac{n}{200}-\b n-3k^2\b n>0$ for $k$ sufficiently large.
\proofend

\section{$K$ has a $k$-factor}

We focus on the subgraph $K$ which we obtained in the previous section by running STRIP on the $k$-core of our random graph, and then possibly deleting one vertex in Section~\ref{sk4}.  We will apply Lemma~\ref{lem:minimality} with $\G:=K$ to prove that $K$ has a $k$-factor.  First we establish some more random graph properties.

\subsection{Further random graph properties}
Recall our setting: we begin with the random graph $G=G_{n,p=c/n}$ where $ c_k+k^{10}\b  \le  c \le c_k+k^{-1/2}$ (see~(\ref{ebeta}),~(\ref{eq:conditions_for_c})).  $K$ is the subgraph of the $k$-core of $G$ obtained after applying \textbf{STRIP}, and possibly deleting one more vertex.  \mmmm{We have shown that \whp\ $K$ satisfies properties~(K1--4) (see Section~\ref{sk4}).  We will require the following additional properties:} 

\begin{lemma}\label{lem:properties}  There exist constants $\g,\e_0 >0,k_0 {\in \mathbb{N}}$ such that for any $k\geq k_0$, \whp\ $K$ satisfies:
\begin{itemize}
\item[(P1)]  For every $Y \subseteq V(K)$ with $|Y| \le 10 \epsilon_0 n$, $e(Y) < \frac{ k |Y|}{6000}$.
\item[(P2)]  For every $Y \subseteq V(K)$ with $|Y|\le \frac12 \mm{|K|}$, $e(Y, V(K) \setminus Y) \ge \gamma k |Y|.$
\item[(P3)]  For every disjoint pair of sets $X,Y \subseteq V(K)$ with $|X| \ge \frac{1}{200}|Y|$ and $|Y| \le \epsilon_0 n$, $e(X,Y) < \frac12 \gamma k |X|$.
\item[(P4)]  For every disjoint pair of sets $X, Y \subseteq V(K)$ with $|X|+|Y| \le \epsilon_0 n$, $e(X,Y) < \left( 1 +\frac{1}{2000}\right) |N(X) \cap Y|+ \frac{ k}{100}|X|$.
\item[(P5)]  For every disjoint pair of sets $S,T \subseteq V(K)$ with $|T|< \frac{1}{10}\epsilon_0 n$ and $|S| > \frac{9}{10}\epsilon_0 n$, $e(S,T) < \frac34 k |S|.$
\item[(P6)]  For every disjoint pair of sets $S,T \subseteq V(K)$ with $|T| \ge \frac{1}{10}\epsilon_0 n$, we have $e(S,T) \le k|S|+\frac{3}{4}\sqrt{k \log k}|T|$ and $\sum_{v \in T}d(v) > (k+\frac{7}{8}\sqrt{k \log k})|T|$.  
\end{itemize}
\end{lemma}

\mmmm{
{\bf Remark:} Note that $\g,\e_0$ do not depend on $k$.

These properties all correspond to very similar properties in~\cite{Mike,Pawel}, with the exception of {(P1)} which is very standard in random graph theory. There are no new ideas here, and the proof is a bit lengthy.  So we defer the proof  of Lemma~\ref{lem:properties}  to Section~\ref{sec:properties}.
}

\subsection{Verifying Tutte's condition}
We now assemble our pieces to show that $K$ has a $k$-factor. 
Recall from Definition~\ref{dlh} that $L,H$ are the vertices of degree $k$ and  at least $k+1$, respectively, and that, eg., $T_L, T_H$ denote $T\cap L, T\cap H$, respectively.

As we proved in Section~\ref{sfts}, \whp\ the subgraph $K$ has the following properties.  Note that (K2) is rephrased using the notation of Definition~\ref{dlh}.

\begin{itemize}
\item [(K1)] for every vertex $v\in K$,  $k \le d_K(v) \le 2 k$;
\item [(K2)] every vertex $v\in H$ has at most $\frac{9}{10}k$ neighbours in $L$;
\item [(K3)]  $|K|\geq \frac{n}{3}$;
\item [(K4)]  $k|K|$ is even.
\end{itemize}

\mmmm{In the previous section, we showed that \whp\ $K$ satisfies properties (P1--6), stated above.}

\dmm{Recall also Tutte's condition~\eqref{etutte2}: 
A graph $\G$ with minimum degree at least $k\geq 1$ has a $k$-factor if \mm{and only if} for every pair of disjoint sets $S,T \subseteq V(\G)$,
$$
k|S|+\sum_{v\in T_H}(d_{\G}(v)-k) \ge q(S,T)+e(S,T),  
$$
and recall also \mm{from Lemma~\ref{lem:minimality}} that a graph $\G$ with minimum degree $k\geq 1$ has a $k$-factor if \mm{and only if} this inequality~\eqref{etutte2} holds for every pair of disjoint sets $S,T \subseteq V(\G)$
satisfying:
\begin{itemize}
\item [(M1)] $S \subseteq H$; and
\item [(M2)] every component $Q$ counted by $q(S,T)$ satisfies $Q_H\neq\emptyset$.
\end{itemize}
}

 By our previous lemmas, it will suffice to prove:
\begin{lemma}\label{lkfactor}  If $K$ satisfies properties (K1--4,P1--{6}) then for every pair of disjoint sets $S,T \subseteq V(K)$ satisfying (M1--2) 
inequality~(\ref{etutte2}) holds.
\end{lemma}

In all but one case, we actually prove~(\ref{etutte}); it will be useful to restate it:
\[k|S|+|T_H| \ge q(S,T)+e(S,T),\]
where  $e(S,T)$ is the number of edges  from $S$ to $T$ and $q(S,T)$ is the number of components $Q$ of $K \setminus (S\cup T)$ such that $k|Q|$ and $e(Q,T)$ have different parity. 
Recall that~(\ref{etutte}) implies~(\ref{etutte2})  since every vertex in $T_H$ has degree at least $k+1$.

{\proofstart}
Let $\epsilon_0,\d >0$ be the constants implied by Lemma~\ref{lem:properties}.   {Recall that $\e_0,\d$ are independent of $k$ (see Remark 1 following the statement of Lemma~\ref{lem:properties}).  So when we lower bound $k$ in what follows, our lower bounds can be in terms of  $\e_0,\d$.  We will assume that $k\geq k_0$ from  Lemma~\ref{lem:properties} and so we can assume that properties (P1--6) all hold.}

We will consider two cases depending on the size of $S \cup T$.

\bigskip

\textit{Case 1: $|S|+|T| \le \epsilon_0 n$:}\\
Recall that $q(S,T)$ \mm{is} the number of connected components $Q$ of $K \setminus (S \cup T)$ such that $k|Q|$ and $e(Q,T)$ have different parities.
Denote by $X$ the union of all vertices belonging to connected components $Q $ of $K \setminus (S \cup T)$ that contribute to $q(S,T)$, other than a largest component of $K \setminus (S \cup T)$ (this largest component might or might not contribute to $q(S,T)$, but is neglected in any case; {if more than one component is largest, we pick one of them arbitrarily}).  As $K \setminus (S \cup T)$ has at most one component of size at least $\hf |K|$, we apply {(P2)} setting $Y$ to be any component in $X$. \dm{N}oting that all edges from $Y$ to $V(K)\bk Y$ are edges from $X$ to $S\cup T$ \dm{and s}umming over all such components yields:
\begin{equation}\label{exp3} 
e(X, \mm{S\cup T}) \ge \gamma k |X|.
\end{equation}
If $|X| > \inv{200}(|S|+|T|)$ then we can apply {(P3)} with $(X,Y):=(X,S \cup T)$ to obtain $e(X,S\cup T)<\hf\g k |X|$ which contradicts~(\ref{exp3}).  So we have:
\begin{equation}\label{ex200} 
|X| \leq \inv{200}(|S|+|T|).
\end{equation}

Recalling that $S_H=S$ by (M1), we now turn our attention to the \mm{vertices of $H$}.
(\ref{ex200})~and the fact that we are in Case~1 imply that $|S|+|T_H|+|X_H|\leq |S|+|T|+|X|<2\e_0 n$ and so we can apply {(P1)} to obtain:
\begin{equation}\label{ep0}
e(S\cup T_H \cup X_H)< \frac{k}{6000} |S\cup T_H \cup X_H|.
\end{equation}
\mm{(K1-2)} imply that each vertex in $X_H$ has at least $\frac{k}{10}$ neighbours in $H$.  Every such neighbour must be in $X\cup S\cup T$ and so $e(X_H,S\cup T_H)\geq \frac{k}{10}|X_H|-2e(X_H)$.  This yields 
\[\mm{e(S\cup T_H \cup X_H)\geq e(X_H)+e(S,T_H)+ e(X_H,S\cup T_H)\geq e(S,T_H)+ \frac{k}{10}|X_H|-e(X_H)},\] 
which combined with~(\ref{ep0}) and {(P1)} applied to $X:=X_H$ gives:
\[\frac{k}{6000} (|S|+ |T_H| +| X_H|) \geq e(S,T_H)+ \left(\frac{k}{10}-\frac{k}{6000}\right)|X_H|\geq \frac{599k}{6000}|X_H|.\]
Rearranging allows us to replace~(\ref{ex200}) with a bound only involving the high vertices:
\begin{equation}\label{exh} 
|X_H| \leq \frac{1}{400}(|S|+|T_H|).
\end{equation}

(M2) implies that every component counted by $q(S,T)$, except possibly the one excluded from $X$, contains a vertex of $X_H$. So $q(S,T)\leq |X_H| +1$.  This allows us to bound the RHS of~(\ref{etutte}) as: 
\begin{align}
\nonumber q(S,&T)+e(S,T)\\
\nonumber &\leq |X_H|+e(S,T_H) +e(S,T_L) +1\\
\nonumber &\leq  \frac{1}{400}(|S|+|T_H|) + \left(1+\frac{1}{2000}\right)|N(S)\cap T_H| +\frac{k}{100}|S| +\frac{9k}{10}|S| +1\\
\nonumber & \qquad\qquad\qquad\qquad\qquad\mbox{ by (\ref{exh}), {(P4)}, (M1) and (K2)}\\
&\le \frac{1}{400}|T_H| + \left(1+\frac{1}{2000}\right)|N(S)\cap T_H| +\frac{92k}{100}|S| +1
 \label{esth}
\end{align}

We now split this case into 3 subcases.

\bigskip

\noindent \textit{Case 1a: $S=\emptyset$.} \\
In this case our goal is to show $|T_H|\geq q(S,T)$. \dm{By~\eqref{esth}, since $e(S,T)=0, |S|=0, |N(S) \cap T_H|=0$, we have $q(S,T)\leq 1 + \frac{1}{400}|T_H|$.}
 This yields  $|T_H|\geq q(S,T)$ if $|T_H|\geq 1$.

Thus we can assume $|T_H|=0$ \mm{and so  $q(S,T)\leq 1$, i.e.\ there is at most one component $Q$ in $K\bk (S\cup T)$ which is counted by $q(S,T)$.}
Since $S=\emptyset$, $e(T,Q)=e(T, K\bk T)$ which has the same parity as ${2e(T)+e(T,Q)=}\sum_{v\in T}d_K(v)=k|T|$ since $T_H=\emptyset$.  
(K4) implies that $k|T|$ has the same parity as $k|Q|$. Therefore $k|Q|$ has the same parity as  {$e(Q,T)$ and so} $q(S,T)=0$.  Therefore the LHS and RHS of~(\ref{etutte}) are both 0 and so {the desired inequality} holds.

\bigskip

\noindent\textit{Case 1b: $|T_H|\leq 3|N(S)\cap T_H|$ and  $S\neq\emptyset$.}\\
\noindent(\ref{esth}) and $|N(S)\cap T_H|\leq |T_H|$ imply
\begin{eqnarray}
\nonumber q(S,T)+e(S,T)
&\leq& 1+ k|S| + |T_H| - \left( \frac{k}{20}|S| - \frac{1}{300}|T_H| \right) -\frac{k}{40}|S|\\
&\leq&  k|S| + |T_H| -   \left( \frac{k}{20}|S| - \frac{1}{300}|T_H| \right)
\label{ec1b}
\end{eqnarray}
for $k\geq40$ since $S\neq\emptyset$.
(K1) implies that $|N(S)\cap T_H|\leq |N(S)| \leq 2 k|S|$.  So since we are in Case~1b:
\[\frac{1}{300}|T_H| \leq\inv{100}|N(S)\cap T_H| \leq  \frac{k}{50}|S| \leq \frac{k}{20}|S|.\]
This and~(\ref{ec1b}) imply~(\ref{etutte}).

\bigskip

\noindent\textit{Case 1c: $|T_H| > 3|N(S)\cap T_H|$ and $S\neq\emptyset$.}\\
\noindent(\ref{esth}) implies
\begin{eqnarray*}
q(S,T)+e(S,T)
&\leq& 1+ \frac{19k}{20}|S| + \frac{1}{400}|T_H| + 2|N_S\cap T_H|\\
&\leq&k|S|+|T_H|
\end{eqnarray*}
for $k\geq20$ since $S\neq\emptyset$. This is~(\ref{etutte}).

\bigskip

\textit{Case 2: $|S|+|T| \ge \epsilon_0 n: $} \\
This case follows  the arguments of Cases 3 and 4 of~\cite{Pawel}. We reproduce them here:

If $|T|<\inv{10}\e_0 n$ then $|S|>\frac{9}{10}\e_0 n$ and so $e(S,T)< \frac{3}{4}k|S|$ by~{(P5)}.  We also have $q(S,T)<n<\inv{4}k|S|$ for $k>\frac{40}{9\e_0}$, and this yields~(\ref{etutte}).  

If $|T|\geq\inv{10}\e_0 n$ then we use the bound $q(S,T)<n<\frac{1}{16}\sqrt{k\log k} |T|$ for $k> 25600 /\e_0^2$, and then the two parts of~{(P6)} combine to give~(\ref{etutte2}).  This is the only case where we prove~(\ref{etutte2}) directly rather than~(\ref{etutte}).
{\proofend}

Our main theorem follows immediately:

{\bf Proof of Theorem~\ref{mt}} We prove in Section~\ref{sfts} that, if $c_k+ {k^{10} \b} \leq c \leq c_k+k^{-1/2}$ then \mm{\whp} the subgraph $K$ obtained by STRIP (and possibly deleting one additional vertex) satisfies properties \mm{(K1-4)}.  Lemma~\ref{lem:properties} establishes that \mm{\whp} $K$ satisfies properties \mm{(P1-6)}.  So Lemma~\ref{lkfactor} implies that \mm{\whp} every pair of disjoint vertex sets $S,T\subseteq V(K)$ satisfying \mm{(M1-2)} also satisfy~(\ref{etutte2}).  Lemma~\ref{lem:minimality} now establishes that $K$ has a $k$-factor and so Theorem~\ref{mt} holds for  $c_k+ {k^{10} \b} \leq c \leq c_k+k^{-1/2}$. Containing a $k$-regular subgraph is a monotone {increasing} property and  {$c_k+e^{-k/300} \ge c_k+k^{10}\b$ for $k$ sufficiently large},  so Theorem~\ref{mt} holds for  all $c\geq c_k+e^{-k/300}$.
\proofend

\mmm{
\begin{remark}\label{rokn}
The proof of Lemma~\ref{lw0}(a), below, shows that the size of the $k$-core in fact is \mmmm{$(1-o_k(1))n$}.
So after removing at most $\b n=e^{-k/200} n$ vertices (Lemma~\ref{lstop}),  the $k$-regular subgraph we obtain is of size \mmmm{$(1-o_k(1))n$.}
\end{remark}
}

\section{Proof of Lemma~\ref{lw0}}\label{sec:lw0}

\mmmm{To complete the paper, all that remains is two deferred proofs. First, we present the proof of Lemma~\ref{lw0}, where we establish some straightforward properties of the $k$-core.  Recall the statement:

{\bf Setup for Lemma~\ref{lw0}:}  {\em $k$ is a sufficiently large constant, and  $c_k < c \le \cmax = c_k + k^{-1/2}$.  $C$~is the $k$-core of $G_{n,p=c/n}$.  $\La$ is a uniform configuration with the same degree sequence as $C$.   {Finally,} $W_0,R$ are as defined in Section~\ref{stsp}. \dmmm{For the convenience of the reader, we state the lemma again.}

{\bf Lemma~\ref{lw0}:}  {W.h.p.\ }before the first iteration of STRIP:
\begin{enumerate}
\item[(a)] $|\La|> {0.99} n$;
\item[(b)] ${ 0.99} \frac{n}{k} < |W_0| < {1.01} \frac{n}{k}$;
\item[(c)] the total degree of {the set of} vertices with  degree greater than $2k$ is at most {$e^{-k/6}n$};
\item[(d)] there are at least $\frac{n}{{5}k}$ edges with both endpoints in $W_0$;
\item[(e)] there are at least ${\frac {1}{2}} n$ edges from $W_0$ to $R$;
\item[(f)] there are at least $\inv{3}k n$ edges with both endpoints in $R$;
\item[(g)] $C$ has at most {$e^{-k/3}n$} vertices of degree at most $2k$ and with at least $\hf k$ neighbours in $W_0$;
\item[(h)] at least $\frac{n}{{200}}$ vertices in $R$ have no neighbours in $W_0$.
\end{enumerate}
}
}

\proofstart
First recall from Section~\ref{stsp} that, before the first iteration of STRIP, $W_0$ is the set of vertices with degree $k$ in $\La$, and $R$ is the set of vertices with degree greater than $k$ in $\La$.

Part (a) is well-known; the size of the $k$-core approaches $n$ as $k$ grows (see, for example, the results of Molloy~\cite{Molloy_cores} or Gao~\cite{Gao}). Indeed, in~\cite{Molloy_cores} it is proved that \whp\ the $k$-core has size ${\z} n + o(n)$, \mm{where}
\[ \z = \z(c) = 1 - e^{-x} \sum_{i=0}^{k-1} \frac {x^i}{i!} = e^{-x} \sum_{i\ge k} \frac {x^i}{i!} = \Pr {( \Po(x) \ge k )}, \]
where $x=x(c)$ is the greatest solution to
\begin{equation}\label{eq:c} 
c = f(x) := \frac{x}{1-e^{-x}\sum_{i=0}^{k-2} x^i/i!}. 
\end{equation}
Recall from~(\ref{eck}) that $c_k$ is the minimum value of $f$ over all $x>0$.  Simple analysis of $f$ shows that  there is exactly one value of $x$ for which $f(x)=c_k$; we denote that value by $x_k$.  \mm{Moreover}, for every $c>c_k$  there are exactly two solutions for $x$. It is straightforward to verify that $f'(x) \ge 0$ for $x \ge x_k$ and so for $c\geq c_k$ we have $x_k\le x(c) \le x(\cmax)$.

{Recall from~(\ref{eq:ck}) that $q_k=\log k-\log (2\pi)$. \cite{Pawel} shows} that 
\begin{equation}\label{eq:xk}
x_k =k+(kq_k)^{1/2}+\frac{q_k}{3}-1+ o_k(1),
\end{equation}
(note that $x_k$ is denoted as $\lambda_k$ in~\cite{Pawel}). It follows that 
$$
\z = \Pr{ ( \Po(x) \ge k )} \ge \Pr {( \Po(x_k) \ge k )} \to 1,
$$
as $k \to \infty$. Part (a) holds for $k$ large enough.

\medskip

Corollary 3 of~\cite{cainwor} establishes that for any constant $i\geq k$, the number of vertices of degree $i$ in the $k$-core is \whp\  $\la_i n +o(n)$ where
\begin{equation}\label{eq:deg_distribution}
\la_i=  \Pr {( \Po(x) = i )} = \frac {e^{-x} x^i}{i!}.
\end{equation}
In particular, \whp\ $|W_0| =  \frac{e^{-x}x^{k}}{k!}n +o(n)$, and so to prove part (b) we will estimate $\frac{e^{-x}x^{k}}{k!}$. 

Below, setting $\d=\d(c)=x(c)-x_k$, we will prove that $c\leq \cmax$ implies
\begin{equation}\label{edbound}
\d\leq \log k.
\end{equation}
So for now, we will restrict our attention to  $x = x_k + \delta$ with $0 \le \delta \le \log k$. Using $1+y = \exp \left( y + O(y^2) \right)$, and  $1/(1+y) = 1 - y + O(y^2)$ we get that
\begin{align}
\frac{e^{-x}x^{k}}{k!} &= \frac{e^{-x_k-\delta}x_k^{k}}{k!} \ \left( 1 + \frac {\delta}{x_k} \right)^k = \frac{e^{-x_k}x_k^{k}}{k!} \ \exp \left( -\delta + \frac {\delta k}{x_k} + O \left( \frac {\delta^2 k}{x_k^2} \right) \right) \nonumber \\
&= \frac{e^{-x_k}x_k^{k}}{k!} \ \exp \left( -\delta + \delta \left( 1 - (q_k/k)^{1/2} + O(q_k/k) + O(\delta / k) \right) \right) \qquad\qquad \mbox{ by~(\ref{eq:xk}) }\nonumber \\
&= \frac{e^{-x_k}x_k^{k}}{k!} \ \exp \left( - \delta (q_k/k)^{1/2} + O(\delta \log k / k) \right)\qquad\qquad\qquad\qquad \mbox{ since } q_k=\log k +O(1)\nonumber \\
&= \frac{e^{-x_k}x_k^{k}}{k!} \left( 1 - \delta (q_k/k)^{1/2} + O(\log^3 k / k) \right) \label{eq:error1} \\
& = \frac{e^{-x_k}x_k^{k}}{k!} (1+o_k(1)), \qquad\qquad\qquad\qquad\qquad\qquad\mbox{ for }\d\leq\log k.\nonumber
\end{align} 
Using Stirling's formula,~(\ref{eq:xk}) and the fact that $1+y=\exp \left( y-y^2/2+y^3/3+O(y^4) \right)$, we get for some $\epsilon = \epsilon(k) = o_k(1)$
\begin{align}
\frac{e^{-x_k}x_k^{k}}{k!} &= \frac{e^{-x_k}}{\sqrt{2\pi k}}\left(\frac{ex_k}{k}\right)^{k} \left( 1+O(k^{-1})\right) \nonumber\\
&= \frac{e^{-(kq_k)^{1/2}-\frac{q_k}{3}+1-\epsilon}}{\sqrt{2\pi k}}\left(1+\frac{(kq_k)^{1/2}+\frac{q_k}{3}-1+\epsilon}{k} \right)^{k} (1+O(k^{-1})) \nonumber \\ 
&= \frac{ \exp \left( - \frac{q_k}{2} + (1-\e)(q_k/k)^{1/2} + O( q_k^2 / k ) \right) }{\sqrt{2\pi k}} \nonumber\\
& =\frac{1+ (1-\e) (q_k/k)^{1/2} + O( \log^2 k / k ) }{k} \ \ 
\mbox{ since } q_k=\log(k/2\pi) \mbox{ and } e^y=1+y+O(y^2) \label{eq:error2} \\
&=\frac{1+ o_k(1)}{k}. \nonumber 
\end{align}
Therefore
\begin{equation}\label{eok1}
\frac{e^{-x}x^{k}}{k!}=\frac{1+ o_k(1)}{k}\qquad\qquad\mbox{ for $x=x_k+\d$ with } 0\leq\d\leq\log k.
\end{equation}
Next we prove that our upper bound $c \le \cmax = c_k + k^{-1/2}$ implies~(\ref{edbound}) and so~(\ref{eok1}) establishes part (b).  To do this we  estimate the derivative of $f(x)$ over the range $0\leq \d\leq \log k$ to show that at $\d = \log k$ and $x=x_k+\d$ we have $c(x)=f(x)>\cmax$.  Rewriting~(\ref{eq:c}) as: 
\[ f(x) = \frac{x}{e^{-x}\sum_{i\geq k-1}\frac{x^i}{i!}},\]
the derivative is
\begin{eqnarray}
f'(x) &=& \frac { e^{-x} \sum_{i \ge k-1} \frac {x^i}{i!} - x \left(-e^{-x} \sum_{i \ge k-1} \frac {x^i}{i!} +e^{-x} \sum_{i \ge k-1} \frac {x^{i-1}}{(i-1)!}\right)} { \left( e^{-x} \sum_{i \ge k-1} \frac {x^i}{i!} \right)^2} \nonumber \\
&=& \frac { e^{-x} \sum_{i \ge k-1} \frac {x^i}{i!} - x e^{-x} \frac {x^{k-2}}{(k-2)!} } { \left( e^{-x} \sum_{i \ge k-1} \frac {x^i}{i!} \right)^2}. \label{ed0}
\end{eqnarray}
Recalling $q_k=\log k-\log(2\pi)$, \mm{and that $x=x_k+\d$ with} $0\leq \d\leq \log k$, \mm{we have}
\begin{align}
x e^{-x} \frac {x^{k-2}}{(k-2)!} &= e^{-x} \frac {x^{k}}{k!} \cdot \frac {k(k-1)}{x} 
= \frac{e^{-x_k}x_k^{k}}{k!} \left( 1 - \delta (q_k/k)^{1/2} + O(\log^3 k / k) \right) \frac {k(k-1)}{x} 
\qquad\mbox{ by~(\ref{eq:error1})} \nonumber \\
&= \left( 1 + (1- \e - \delta) (q_k/k)^{1/2} + O(\log^3 k/k) \right) \frac {k-1}{x}
\qquad \qquad\qquad\qquad\qquad\qquad \mbox{ by~(\ref{eq:error2})} \nonumber \\
&= \left( 1 + (1- \e - \delta) (q_k/k)^{1/2} + O(\log^3 k/k) \right) 
\left( 1 -(q_k/k)^{1/2} + O(\mm{\log k}/k) \right) \nonumber \\
 &\qquad\qquad\qquad\qquad\qquad\qquad\qquad\qquad \mbox{ by~(\ref{eq:xk}) and since } \inv{1+y} = 1 - y + O(y^2) \nonumber \\
&= 1 - (\e +\delta) (q_k/k)^{1/2} + O(\log^3 k/k).\label{ed7}
\end{align}
Standard bounds  for the tail probabilities of a Poisson random variable  (see eg.~\cite[Theorem A.1.15]{AS}), 
along with~(\ref{eok1}) yield that for $0\leq \d\leq \log k$ we have
\[
\Pr{ \left( \Po(x) \le k-2 \right) } 
= O \left( k^{-1/2} \right),
\]
and so
\begin{equation}\label{ed2}
e^{-x} \sum_{i \ge k-1} \frac {x^i}{i!} = 1 - \Pr{ \left( \Po(x) \le k-2 \right) } = 1 - O \left( k^{-1/2} \right).
\end{equation}
Substituting~(\ref{ed7}) and~(\ref{ed2}) \dm{into~(\ref{ed0})} and recalling $q_k=\log k +O(1)$ yields
\begin{eqnarray}
f'(x) &=& (\e +\delta) (q_k/k)^{1/2} + O \left( k^{-1/2} \right) \nonumber \\
&=& (\e +\delta) (\log k/k)^{1/2} + O \left( k^{-1/2} \right),
\mbox{ for $x=x_k+\d$ and } 0\leq\d\leq\log k. \label{eq:bound_for_f'}
\end{eqnarray}
Therefore, recalling that $f(x_k)=c_k$ and since $\e=o_k(1)=o(\d)$ we have for $k$ sufficiently large
\[
f(x_k + \log k) = f(x_k) + \Theta( (\log k)^2 )\cdot (\log k / k)^{1/2} +(\log k)\cdot O \left( k^{-1/2} \right)
> c_k + k^{-1/2}=\cmax.
\]
\dm{(Note that integration of~\eqref{eq:bound_for_f'} would have given the same result.)} Since $f(x)$ is monotone increasing for $x>x_k$, it follows that $x(c)\leq x_k+{\log k}$ for all $c_k\leq c\leq \cmax$, thus proving~(\ref{edbound}). Therefore~(\ref{eok1}) and~(\ref{eq:deg_distribution}) yield part (b).

\medskip

For part (c) it is easier to show the desired property for $G$ instead of $C$; the conclusion for $C$ will trivially follow. The degree of each vertex in $G$ is a random variable $X$ with the binomial distribution $\Bin(n-1,c/n)$ with $\E[X] \le c < 1.1 k$, provided that $k$ is large enough. Hence, \mmm{applying~\eqref{chernoff:up} with $t=0.9k+\ell$}, we get that the expected total degree of all vertices with degree greater than $2k$ is at most 
\begin{align*}
n \sum_{\ell \ge 1} (2k + \ell) \Pr{ ( X \ge 2k + \ell )} & \le n \sum_{\ell \ge 1} (2k+\ell) \exp \left( - \frac {(0.9k+\ell)^2}{\mmm{2 (1.1k+ (0.9k+\ell)/3)} } \right) \\
& \le n \sum_{\ell \ge 1} (2k+\ell) \exp \left( - \frac {k+\ell}{4} \right) = O( k e^{-k/4} n ) < e^{-k/5} n,
\end{align*}
provided $k$ is sufficiently large. Since the concentration can be proved with a straightforward concentration argument using, for example, Azuma's Inequality or an easy second moment argument, we omit the details. 
This establishes part (c).

\medskip

For the remaining parts, let us first observe that the total degree of all vertices in $G$ is $2 \Bin \left( {n \choose 2}, c/n \right)$ with expectation $c(n-1) < 1.01 kn$, provided $k$ is sufficiently large. Hence, by Chernoff's bound, \whp\ it is at most $1.02 kn$, and this upper bound clearly holds for the total degree of the vertices of $C$. On the other hand, part (a) implies that it is at least $0.99 kn$. 

For parts (d,e,f), \mm{we return to analyzing $\La$ directly. We will} focus on the partners of the vertex-copies in $W_0$.
Part~(b) implies that the total degree of the vertices in $W_0$ is between $0.99n$ and $1.01n$.  Expose the partners in $\La$ of the vertex-copies of $W_0$, one at a time.  At step $i\leq 1.01 n$, the probability that the partner chosen is in $W_0$ is between $p_i$ and $q_i$, where
\begin{equation}\label{epiqi}
p_i = \max \left\{ \frac {0.99n - 2i}{1.02kn}, 0 \right\} 
\text{ and } 
q_i = \max \left\{ \frac {1.01n - i}{0.99kn-2i}, 0 \right\} \le \max \left\{ \frac {1.01n - i}{0.98kn}, 0 \right\},
\end{equation}
provided $k$ is sufficiently large. Hence, the number of vertex-copies in $W_0$ whose partner is in $W_0$ can be stochastically lower/upper bounded by the two sums of independent Bernoulli random variables: the first one with parameters $p_i$ and the second one with $q_i$. The expected value of the first sum is more than $0.24 n / k$, and the expectation of the second one is at less than $0.53 n / k$. The concentration follows immediately from Chernoff's bound and we get that \whp\ the number of edges with both endpoints in $W_0$ is at least $\frac {n}{5k}$ and at most $\frac {n}{k}$, which finishes part (d).

\medskip

Now, parts~(e) and~(f) follow deterministically. The number of edges from $W_0$ to $C\bk W_0$ is at least $0.99 n - 2 \frac {n}{k} \ge \frac {1}{2} n$ for $k$ large enough. Finally, there are at least $(0.99 kn - 2 \cdot 1.01 n)/2$ edges with both endpoints in $C\bk W_0$ which is more than $\frac {1}{3} kn$ for $k$ large enough.

\medskip

Part~(g) is slightly more complicated. For a contradiction, suppose that there \dm{is a set $T$ with $|T|=e^{-k/3}n$ 
such that every vertex in $T$ is} of degree at most $2k$ and \dm{has} at least $\hf k$ neighbours in $W_0$. (Note that some of them might be from $W_0$.) \dmm{We will first show that t}his implies (deterministically) that there exists a \mmm{subset $S\subseteq T$ of size at least $\hf e^{-k/3}n$ and} with at least $\frac {1}{8} k e^{-k/3}n$ edges between $S$ and $W_0 \setminus S$. 


\mmmm{To prove this, we consider the average of $|E(S,W_0\setminus S)|$ over all subsets $S\subset T$ of size $\rup{\hf|T|}$.   Consider any edge $uv$ with $u\in T$ and $v\in W_0$; there are at least $\inv{4}k|T|$ such edges (if $u,v$ are both in $T\cap W_0$ then the edge $uv$ is only counted once). A very simple count shows that $uv\in E(S,W_0\setminus S)$ for at least half of the subsets $S$ of size $\rup{\hf|T|}$; there are four cases corresponding to the parity of $|T|$ and whether $v\in T$. Therefore,
\[\sum_{S\subset T, |S|=\rup{\hf|T|}}|E(S,W_0\setminus S)|\geq\inv{8}k|T|{|T|\choose \rup{\hf|T|}},\]
and so at least one such set $S$ has $|E(S,W_0\setminus S)|\geq\inv{8}k|T|$.

Next, we show that \whp\ no such set $S$ exists in $\Lambda$.}
Indeed, let us fix \dm{a} set $S$ and expose \mm{the} partners of all vertex-copies from $W_0 \setminus S$ in $\La$. Arguing similarly as for~(\ref{epiqi}), the number of edges from $W_0 \setminus S$ to $S$ can be stochastically upper bounded by the binomial random variable
$$
X \sim \Bin \left( 1.01 n, \frac {2k|S|}{0.98 kn} \right), \quad \text{ with } \quad \E[X] < 3 e^{-k/3} n.
$$
It follows from Chernoff's bound that the expected number of sets $S$ with a large number of edges to $W_0 \setminus S$ is at most
\begin{align*}
{n \choose \hf \dmmm{\lceil}e^{-k/3}n\dmmm{\rceil}} & \Pr{ \left( X \ge \frac {ke^{-k/3}}{8} n \right) } \le \left( \dmmm{3}e^{1+k/3} \right)^{\hf \dmmm{\lceil}e^{-k/3}n\dmmm{\rceil}} \exp \left( - 1.4 \cdot \frac {k e^{-k/3}}{8} n \right) = o(1),
\end{align*}
provided that $k$ is sufficiently large. Part~(g) follows from Markov's inequality.

\medskip

Finally, let us move to part~(h). We showed earlier that the total degree of the vertices of $C$ is at least $0.99kn$ and the total degree of the vertices in $W_0$ is at most $1.01n$. It follows from parts~(a), (b), and~(c) that there are at least $0.98n$ vertices in $R$ that are of degree at most $2k$ (and, of course, at least $k+1$). We pick (arbitrarily) $0.13n$ of them and expose partners of all corresponding vertex-copies. Note that, regardless of the history of the process, the probability that a given vertex of degree $\ell \le 2k$ has no neighbour in $W_0$ is at least
\[ 
\left( 1 - \frac {1.01n}{0.99 kn - (2k)(0.13n)} \right)^{\ell} \ge \left( 1 - \frac {1.4}{k} \right)^{2k} \ge e^{-3},
\]
provided that $k$ is sufficiently large. Hence, the number of vertices in $R$ that have no neighbours in $W_0$ is bounded from below by the random variable $X \sim \Bin(0.13 n, e^{-3})$ with $\E[X] = 0.13 e^{-3} n > 0.006 n$. Part (h) holds by Chernoff's bound and the proof of the lemma is finished.
\proofend

\section{Proof of Lemma~\ref{lem:properties}}\label{sec:properties}

\mmmm{
Our final piece is the deferred proof regarding properties of the subgraph $K$ obtained by our stripping procedure.
Recall our setting: we begin with the random graph $G=G_{n,p=c/n}$ where $ c_k+k^{10}\b  \le  c \le c_k+k^{-1/2}$ (see~(\ref{ebeta}),~(\ref{eq:conditions_for_c})).  $K$ is the subgraph of the $k$-core of $G$ obtained after applying \textbf{STRIP}, and possibly deleting one more vertex. We repeat the statement of the lemma:

{\bf Lemma~\ref{lem:properties}} There exist constants $\g,\e_0 >0,k_0 {\in \mathbb{N}}$ such that for any $k\geq k_0$, \whp\ $K$ satisfies:
\begin{itemize}
\item[(P1)]  For every $Y \subseteq V(K)$ with $|Y| \le 10 \epsilon_0 n$, $e(Y) < \frac{ k |Y|}{6000}$.
\item[(P2)]  For every $Y \subseteq V(K)$ with $|Y|\le \frac12 |K|$, $e(Y, V(K) \setminus Y) \ge \gamma k |Y|.$
\item[(P3)]  For every disjoint pair of sets $X,Y \subseteq V(K)$ with $|X| \ge \frac{1}{200}|Y|$ and $|Y| \le \epsilon_0 n$, $e(X,Y) < \frac12 \gamma k |X|$.
\item[(P4)]  For every disjoint pair of sets $X, Y \subseteq V(K)$ with $|X|+|Y| \le \epsilon_0 n$, $e(X,Y) < \left( 1 +\frac{1}{2000}\right) |N(X) \cap Y|+ \frac{ k}{100}|X|$.
\item[(P5)]  For every disjoint pair of sets $S,T \subseteq V(K)$ with $|T|< \frac{1}{10}\epsilon_0 n$ and $|S| > \frac{9}{10}\epsilon_0 n$, $e(S,T) < \frac34 k |S|.$
\item[(P6)]  For every disjoint pair of sets $S,T \subseteq V(K)$ with $|T| \ge \frac{1}{10}\epsilon_0 n$, we have $e(S,T) \le k|S|+\frac{3}{4}\sqrt{k \log k}|T|$ and $\sum_{v \in T}d(v) > (k+\frac{7}{8}\sqrt{k \log k})|T|$.  
\end{itemize}

}

{\proofstart}
For every property except {(P2)}, we actually show that it holds in $G$, and so we can work in the $G_{n,p}$ model.
For each property, we will show that it holds if $\e_0$ is sufficiently small.  Thus we can take a value of $\e_0$ that is sufficiently small for all properties.

The following property follows from Lemma 3 in~\cite{Pawel}, where they show that in fact  there is no such subgraph in $G$:

\noindent (P0) For every $Y \subseteq V(K)$ with $|Y| \le 2\log n/(e c \log \log n)$, $e(Y) \le |Y|.$

Now, {(P1)}  follows from a very standard first moment argument applied to $G=G_{n,p=c/n}$, so long as $\e_0$ is sufficiently small.  For $Y \subseteq V(K)$ with $|Y|=s \le 2\log n/(e c \log \log n)$, the statement follows immediately by Property~{(P0)}, for $k\geq 6000$. For larger values of $s$,  and for $k\geq 1f$ and since $k<c<2k$, the expected number of sets $Y\subset V(G)$ with $|Y|=s$ for which $Y$ contains at least $ \frac{k s}{6000}$ edges in $G$ is at most
\begin{eqnarray*}
\binom{n}{s} \binom{\binom{s}{2}}{\frac{ ks}{6000}} p^{\frac{ ks}{6000}} &\le& \left(\frac{ne}{s}\right)^s  \left(\frac{300cse}{ kn}\right)^{\frac{ks}{6000}} \le \left(\frac{ne}{s}\right)^s  \left(\frac{6000se}{ n}\right)^{2s} \\
&=& \left(\frac{6000^2 e^3 s}{ n}\right)^{s} 
< 2^{-s},
\end{eqnarray*}
for $s\leq \e_0n$ so long as $\e_0<1/(2\times6000^2 e^3)$.

Summing over all $2\log n/(e c \log \log n) < s \le 10 \epsilon_0 n$, we obtain that the expected number of sets that fail the desired property is $O(2^{-2\log n/(e c \log \log n)}) = o(1)$. Property {(P1)} now follows by Markov's inequality. 

\smallskip

For property {(P2) we consider two cases}:  

{\em Case 1: $|Y|\leq10\e_0n$.}   Since each vertex in $Y$ has degree at least $k$, and applying {(P1)} we know $e(Y,V(K)\bk Y)\geq k|Y|-2e(Y)\geq k|Y| - k|Y|/3000>\hf k |Y|$.

{\em Case 2:  $|Y|>10\e_0n$.}  Lemma 2 of \cite{Pawel} proves that in the $k$-core $C$, $e(Y,V(C)\bk Y)\geq \g'k|Y|$ for some constant $\g'>0$ independent of $k$. (In fact, they prove this for the $(k+2)$-core but the same proof applies to the $k$-core; the main tool is Lemma 5.3 of {Benjamini, Kozma and Wormald}~\cite{Benjamini}.)  By Lemma~\ref{mtsg},  $|C\bk K|\leq \b n$.  Since $K$ is an induced subgraph of $C$ and every vertex of $Y$ has degree at most $2k$  (by property (K1)), this implies that  $e(Y,V(K)\bk Y)\geq \g'k|Y|-2k\b n$.  This is at least $\hf \g' k|Y|$ if  $k$ is sufficiently large so that $2k\b < 5\g'\e_0$ (recall $\b=e^{-k/200}$).  So {(P2)} holds for $\g=\hf\g'$.

\smallskip

Property {(P3)} follows from Property (P4) of~\cite{Mike}. There it is shown that in $G(n,p)$ with $p=c/n$ and $0  < c  < 2k$ there is no subgraph satisfying the desired property. Their inequality is not strict, but it can be clearly made strict. Their proof holds for any $\g>0$ so long as $\e_0$ is sufficiently small in terms of $\g$.  Thus, we can use the same value of $\g$ as in property {(P2)}.

\smallskip

For property {(P4)},
 it clearly suffices to only consider disjoint sets $X,Y\neq\emptyset$ for which 
 $N(X) \cap Y=Y.$
 \mmm{So it suffices to prove that 
 \whp\ every disjoint pair of sets $X,Y \subseteq V(G)$ with $G \in \mathcal{G}(n,p)$ with $p=c/n$ and $0  < c  < 2k$ and $2\leq |X|+|Y| \le \epsilon_0 n$ satisfies 
\begin{equation}\label{eq:P5}
e(X,Y) \le \mm{\left(1+\frac{1}{2000}\right)}|Y|+ \frac{k}{100}|X|.
\end{equation}
}

The argument is essentially the same as for Lemma 2.5 of~\cite{Mike}, but with different constants; we give it here for the sake of completeness. 

 Let $\s n=|X|$ and $\t n=|Y|$. For any choice of $\sigma, \tau$, the expected number of  sets $X,Y$ in $G$ violating~\eqref{eq:P5} is at most 
 \begin{align}
 \nonumber \binom{n}{\sigma n} & \binom{n}{\tau n}\binom{(\sigma n)(\tau n)}{(1+\frac{1}{2000}) \tau n+ \frac{k}{100}\sigma n}\left(\frac{c}{n}\right)^{(1+\frac{1}{2000} )\tau n+ \frac{ k}{100}\sigma n}  \\
& \le  \left( \frac{e}{\sigma}\right)^{\sigma n}\left( \frac{e}{\tau}\right)^{\tau n}  \left( \frac{e \sigma \tau c}{(1+\frac{1}{2000} )\tau +\frac{ k}{100}\sigma}\right)^{(1+\frac{1}{2000} )\tau n+\frac{ k}{100}\sigma n}.
\label{ep5}
\end{align}
If $\s n, \t n$ are both less than $\sqrt{n}$; i.e. $\s,\t<n^{-1/2}$ then~(\ref{ep5}) is at most
\begin{align}
\nonumber
  \left( \frac{e}{\sigma}\right)^{\sigma n} & \left( \frac{e}{\tau}\right)^{\tau n} 
    \left( \frac{e \sigma \tau c}{\left(1+\frac{1}{2000} \right)\tau} \right)^{\frac{ k}{100}\sigma n}
 \left( \frac{e \sigma \tau c}{\frac{ k}{100}\sigma}\right)^{(1+\frac{1}{2000} )\tau n}\\
 \nonumber
 &=   \left( \frac{e}{\sigma}\right)^{\sigma n}
    \left( \frac{e \sigma c}{1+\frac{1}{2000} } \right)^{\frac{ k}{100}\sigma n}
   \left( \frac{e}{\tau}\right)^{\tau n} 
 \left( \frac{100e  \tau c}{k}\right)^{(1+\frac{1}{2000} )\tau n}\\
 \nonumber & <(A n^{-1/2})^{(\frac{ k}{100}-1)\sigma n} (B n^{-1/2})^{\frac{1}{2000}\tau n}
 \qquad\qquad \mbox{for \mm{some} constants $A,B$}\\
 &<n^{-\frac{1}{4}\frac{1}{2000}(\s n + \t n)}, \qquad\qquad\qquad\qquad\mbox{ for $k>200$.}
 \label{ep5b}
 \end{align}
For general $\s n, \t n$, using $c < 2k$ and $\tau < \epsilon_0$ for $\e_0$ sufficiently small, we see that the \mm{base of the exponent in the} third factor of~(\ref{ep5}) is at most
$$
 \frac{e \sigma \tau c}{(1+\frac{1}{2000}) \tau +\frac{k}{100}\sigma}< \frac{e\sigma \tau c}{\frac{ k}{100}\sigma} < 200e\tau < \left(\frac{\tau}{e^2}\right)^{\frac{1}{ 1+\frac{1}{2000}}}.
 $$
 If $\sigma \ge e^{- k/200}$, then for $k\geq 400$  we have
 $$
  \frac{e \sigma \tau c}{\left(1+\frac{1}{2000} \right)\tau + \frac{k}{100}\sigma}<\left(\frac{\tau}{e^2}\right)^{\frac{1}{ 1+\frac{1}{2000}}}  <e^{-1} \leq\left(\frac{\sigma}{e^2}\right)^{\frac{100}{ k}},
 $$
 while if $\sigma < e^{-k/200}$, then for $k\geq 40000$ we have
 $$
   \frac{e \sigma \tau c}{\left(1+\frac{1}{2000} \right)\tau + \frac{k}{100}\sigma}< \frac{e\sigma \tau c}{\left(1+\frac{1}{2000} \right)\tau} < ec \sigma^{1/2}\sigma^{1/2} < 2ek e^{-k/400}\sigma^{1/2}< \sigma^{1/2} < \left(\frac{\sigma}{e^2} \right)^{\frac{100}{k}}.
 $$
 Hence, the expected number of disjoint sets $X,Y$ with $|X|=\sigma n$, $|Y|=\tau n$ \dm{and $e(X,Y) \ge \left(1+\frac{1}{2000} \right)\mmm{|Y|}+\frac{k}{100}|X|$} is at most 
 \begin{equation}\label{ep5c}
\left( \frac{e}{\sigma}\right)^{\sigma n}\left( \frac{e}{\tau}\right)^{\tau n} 
\left(\frac{\sigma}{e^2} \right)^{\frac{100}{ k} \frac{k}{100}\sigma n}\left( \frac{\tau}{e^2}\right)^{\frac{1}{ 1+\frac{1}{2000}}\left(1 +\frac{1}{2000} \right)\tau n}=\left(\frac{1}{e}\right)^{|X|+|Y|}.
\end{equation}
  For every choice of $y=|X|+|Y|$, there are $y-1$ choices for $|X|, |Y|\geq 1$. Applying~(\ref{ep5b}) and~(\ref{ep5c}), the expected number of sets violating~\mmm{(\ref{eq:P5})} is at most 
  $$
   \sum_{y=2}^{\sqrt{n}} (y-1)n^{-\frac{1}{4}\frac{1}{2000}y} ~~+~~
  \sum_{y=\sqrt{n}}^{\epsilon_0 n} (y-1) \left(\frac{1}{e}\right)^{y}=o(1),
  $$
and, by Markov's inequality,~\mmm{(\ref{eq:P5}) holds and thus Property {(P4)} holds.}
  
 \smallskip 

Property {(P5)}  follows in the same way as property (P8) of~\cite{Mike}, which used proofs taken from~\cite{Pawel}.  The proofs hold for any $\e_0$ sufficiently small so long as $k$ is sufficiently large in terms of $\e_0$.

 \smallskip 

Property {(P6)} comes from Case 4 of the proof of Theorem 1 in~\cite{Pawel}.  The first part is equation (14) of that paper with $\e:=\inv{4}$; it is easy to check that their proof goes through with that value of $\e$. The second part is from the line preceding (14) with $\e:=\inv{8}$; that line holds for every $\e$ so long as $k$ is sufficiently large in terms of $\e$.  In both cases, the proof analyzes the degree sequence of $G_{n,p}$ and so holds for every subset $S,T$ of $G_{n,p}$ so long as the vertices of $T$ have degree at least $k$.  There is a minor difference in that they have $c>c_{k+2}$ rather than $c>c_{k}$ but this has no significant effect on the proof.  Again, the proofs hold for any $\e_0$ sufficiently small so long as $k$ is sufficiently large in terms of $\e_0$.
{\proofend}

\section*{Acknowledgements} Part of this research was conducted while M. Molloy was an Invited Professor at  the \'Ecole Normale Sup\'erieure, Paris and while D. Mitsche was visiting Ryerson University.
The authors are supported by NSERC Discovery grants and an NSERC Engage grant.  \mm{We are grateful to two referees who read this paper very carefully and provided many helpful comments and corrections.}

\end{document}